\documentclass[a4paper]{article}

\usepackage[utf8]{inputenc} \usepackage{lmodern} \usepackage[ngerman, english]{babel} \usepackage{csquotes} 

\usepackage[unicode, colorlinks]{hyperref} \usepackage{enumitem}

\usepackage{emptypage} \usepackage[dvipsnames]{xcolor}

\usepackage{wrapfig} 

\usepackage[ruled,linesnumbered]{algorithm2e}


\usepackage{mathtools} \usepackage{amsthm, amssymb} \usepackage{aligned-overset}

 \usepackage[numbers]{natbib}
\theoremstyle{plain}\newtheorem{prop}{Proposition}[section]
\newtheorem{lemma}[prop]{Lemma}
\newtheorem{corollary}[prop]{Corollary}
\newtheorem{theorem}[prop]{Theorem}

\theoremstyle{definition}
\newtheorem{definition}[prop]{Definition}
\newtheorem{example}[prop]{Example}
\newtheorem{counterexample}[prop]{Counterexample}

\theoremstyle{remark}
\newtheorem{remark}[prop]{Remark}

\makeatletter

\newcommand*{\dims}{d}
\newcommand*{\ind}{\mathbf{1}}

\newcommand*{\inter}{ \mathchoice{\mskip1.5mu{:}\mskip1.5mu}{\mskip1.5mu{:}\mskip1.5mu}{{:}}{{:}}
}

\newcommand*{\real}{\mathbb{R}}
\newcommand*{\nat}{\mathbb{N}}
\newcommand*{\integer}{\mathbb{Z}}

\newcommand*{\domain}{\mathbb{X}}
\newcommand*{\range}{\mathbb{Y}}
\newcommand*{\cX}{\mathcal{X}}
\newcommand*{\cY}{\mathcal{Y}}

\newcommand*{\metric}{d}
\newcommand*{\closure}[1]{\overline{#1}}
\newcommand*{\interior}[1]{\mathrm{int}\, #1}

\newcommand{\E}{\mathbb{E}}
\renewcommand{\Pr}{\mathbb{P}}
\newcommand*{\filt}{\mathcal{F}}
\newcommand*{\borel}{\mathcal{B}}
\newcommand*{\cF}{\mathcal{F}}
\newcommand*{\cG}{\mathcal{G}}
\newcommand*{\cE}{\mathcal{E}}
\newcommand*{\cA}{\mathcal{A}}

\DeclareMathOperator{\Cov}{Cov}
\newcommand*{\normal}{\mathcal{N}}
\newcommand*{\uniform}{\mathcal{U}}

\newcommand*{\indep}{\perp\!\!\!\perp}
\newcommand*{\as}{\text{a.s}\@ifnextchar.{}{\text{.\@} }} \newcommand*{\kernel}{\kappa}

\newcommand*{\rf}{\mathbf{f}}
\newcommand*{\rg}{\mathbf{g}}
\newcommand*{\C}{\mathcal{C}}

\newcommand*{\transpose}{T}

\DeclareMathOperator*{\argmin}{arg\,min}
\DeclareMathOperator*{\argmax}{arg\,max}

\newcommand*{\noise}{\varsigma}

\definecolor{primary}{RGB}{0,48,86}

\colorlet{primary10}{primary!10}
\colorlet{primary25}{primary!25}
\colorlet{primary40}{primary!40}
\colorlet{primary55}{primary!55}
\colorlet{primary70}{primary!70}
\colorlet{primary85}{primary!85}

\definecolor{accent}{RGB}{179,182,185}

\definecolor{magenta}{HTML}{D81B66}
\definecolor{blue}{HTML}{1F76C1}
\definecolor{yellow}{HTML}{FFC107}
\definecolor{teal}{HTML}{00B981}

\newcommand*{\blue}[1]{{ #1}}
\newcommand*{\magenta}[1]{{\color{magenta} #1}}
\newcommand*{\teal}[1]{{\color{teal} #1}}

\makeatother 

\title{Measure Theory of Conditionally Independent Random Function Evaluation}
\author{
	Felix Benning\\
	\texttt{felix.benning@uni-mannheim.de}\\
	University of Mannheim
}

\begin{document}
	\maketitle 

	\begin{abstract}
    In sequential design strategies, common in geostatistics and Bayesian
    optimization, the selection of a new observation point \(X_{n+1}\) of a random
    function \(\rf\) is informed by past data, captured by the filtration
    \(\filt_n=\sigma(\rf(X_0),\dots,\rf(X_n))\). The random nature of \(X_{n+1}\) introduces
    measure-theoretic subtleties in deriving the conditional distribution
    \(\Pr(\rf(X_{n+1})\in A \mid \filt_n)\). Practitioners often resort to a
    heuristic: treating \(X_0,\dots, X_{n+1}\) as fixed parameters within the
    conditional probability calculation. This
    paper investigates the mathematical validity of this widespread practice. We
    construct a counterexample to prove that this approach is, in general,
    incorrect. We also establish our central positive result: for continuous
    Gaussian random functions and their canonical conditional distribution, the
    heuristic is sound. This provides a rigorous justification for a
    foundational technique in Bayesian optimization and spatial statistics. We
    further extend our analysis to include settings with noisy evaluations and
    to cases where \(X_{n+1}\) is not adapted to \(\filt_n\) but is
    conditionally independent of \(\rf\) given the filtration.



\end{abstract}

	\section{Introduction}

Researchers have been confronted with the challenge of optimizing an unknown
function across numerous disciplines. The Bayesian approach addresses this
challenge by modelling the unknown function as a random draw from a prior
distribution over functions.
This probabilistic framework ensures
the conditional distribution of a \emph{random function}\footnote{
    While used synonymously, we avoid the more common term `stochastic process' which
     invokes the notion of a one-dimensional index representing `time' and a filtration
     associated to this time. The domain \(\domain\) is generally
    unordered, e.g.\@ \(\domain=\real^\dims\), and the filtration we consider
     naturally arises from the sequence of evaluations of this random function
     \(\rf\).
 } \(\rf=(\rf(x))_{x\in \domain}\)
given a set of function evaluations \(\rf(x_1), \dots, \rf(x_n)\) is
well-defined. This conditional distribution can then be used to choose the next evaluation
point. This approach to optimization has independently emerged across multiple
research domains, each developing its own terminology for very similar methods.
The most prominent is perhaps \emph{Bayesian optimization}
\citep{kushnerNewMethodLocating1964,jonesEfficientGlobalOptimization1998,frazierBayesianOptimization2018,garnettBayesianOptimization2023},
which is best know in machine learning for its effectiveness at hyperparameter
tuning. However similar techniques were already developed earlier in
\emph{geostatistics} known as ``Kriging''
\citep{krigeStatisticalApproachBasic1951,matheronPrinciplesGeostatistics1963,steinInterpolationSpatialData1999}.

\smallskip

In this paper we highlight and address measure-theoretic challenges
 that must be overcome for a rigorous mathematical treatment of random function
 optimization. In particular, we show how non-rigorous but intuitive claims
 about the conditional distributions of the objective function \(\rf\)
 can be made rigorous. In our main result we formalize
 the heuristic that evaluation locations \(X_n\), which are measurable
 with respect to \(\filt_{n-1} = \sigma(\rf(X_0), \dots, \rf(X_{n-1}))\) 
 (``previsible''), may be treated as if they are deterministic during the
 calculation of the conditional distribution
 \[
     \Pr(\rf(X_n)\in \cdot \mid \rf(X_0), \dots, \rf(X_{n-1})).
 \]
If there is a formula for such a conditional distribution
\textbf{for all} deterministic \(X_k\), then we call this
formula a \emph{joint} kernel (Definition \ref{def: joint kernel}).
If this formula is furthermore valid \textbf{for all} previsible random \(X_k\),
then we call it ``plug-in consistent'' (PIC) (see Definition \ref{def: pic
simple} and \ref{def: pic}). Standard disintegration results can only
provide results for almost all deterministic \(X_k\) with respect to
a single distribution over \(X_k\). One may therefore call our results
``strong disintegration results'' (see Remark \ref{rem: strong disintegration}).
In Counterexample \ref{ex: inconsistent joint
conditional distribution} we show that the PIC property is not guaranteed if the
assumptions of our main results are violated.

\begin{remark}[Measurability of random evaluations]
    \label{rem: measurability of random eval}
    A priori, even the measurability of \(\rf(X)\) -- its very
    existence as a random variable -- is not guaranteed.
However, a measurable evaluation map \(e(f,x) \coloneqq f(x)\) would provide
    this guarantee and
we show that this requirement is
    satisfied in considerable generality: If \(\rf\) is a \emph{continuous} random
    function with locally compact, separable, metrizable domain \(\domain\) and
    Polish codomain \(\range\), then the evaluation map \(e\) is continuous, and
    hence measurable (see Theorem~\ref{thm: topology of continuous functions}). This
    covers essentially all continuous applications of interest (e.g.\
    \(\domain \subseteq \mathbb{R}^d\) and \(\range = \mathbb{R}^n\)). The
    difficulty lies in the fact that the evaluation map 
    is only well known to be continuous with respect to the compact-open
    topology, which is different from the product topology used to
    generate the Borel \(\sigma\)-algebra for random functions (Remark \ref{rem:
    top of pointwise convergence}). 
\end{remark}

\paragraph*{Outline}
In Section \ref{sec: measure theory of random function evaluation} we present the theory, first for a simplified case in Section \ref{sec: dependent evaluation},
    then for the general case in Section \ref{sec: random eval in conditional}.
    For the case where the assumptions of these results are violated
    we present a counterexample in Section \ref{sec: counterexample}.
In Section \ref{sec: applications} we present applications of our results.
Section \ref{sec: proofs} contains the proofs for the results from Section
\ref{sec: measure theory of random function evaluation}.
In Section \ref{sec: topological foundation} we address the issues outlined
in Remark \ref{rem: measurability of random eval}. While these results are likely to be known, we
could not find them anywhere in this generality.

 	\section{Measure theory of random function evaluation}
\label{sec: measure theory of random function evaluation}

Throughout the paper we assume there exists an underlying probability space
\((\Omega, \cA, \Pr)\).

\subsection{Simplified case: Dependent evaluation}
\label{sec: dependent evaluation}

Let \(\rf=(\rf(x))_{x\in \domain}\) be a random continuous function with locally
compact, separable metrizable domain \(\domain\) and Polish codomain \(\range\).
That is we assume \(\rf\) is a random variable in the space of continuous functions
\(C(\domain, \range)\).
In this section we consider the simplified case
\begin{equation}
    \label{eq: measurable evaluation I}    
    \Pr(\rf(X) \in A \mid \filt)
\end{equation}
with only a single random evaluation point \(X\) in \(\domain\) that is measurable
with respect to a sub-\(\sigma\)-algebra \(\filt\) of \(\cA\) and sets \(A\) in \(\borel(\range)\).
Here \(\borel(\range)\) denotes the Borel \(\sigma\)-algebra on \(\range\).
Assume we have access to a formula to calculate \eqref{eq: measurable evaluation I}
for any deterministic \(X\). Specifically, assume we have
a collection \((\kernel_x)_{x\in \domain}\) of regular conditional
distributions such that for all \(x\in\domain\), \(A\in \borel(\range)\) and for \(\Pr\)-almost all \(\omega \in \Omega\),\footnote{
    \label{footnote: evaluation of conditional probability}
    Recall that \(\Pr(\rf(x) \in A\mid \filt)\) is only defined as an \(L^1\)-equivalence class of random
    variables. For an equivalence class \([Z]\) in \(L^1\) it is standard to
    say, ``\([Z](\omega) = Y(\omega)\) holds for almost all \(\omega\)'', if for
    all \(Z\in [Z]\) there exists a null set \(N_Z\) such that
    \(Z(\omega)= Y(\omega)\) for all \(\omega\in N_Z^\complement\). Observe that
    the null set \(N_Z\) depends on \(Z\). For
    \(\Pr(\rf(x)\in A \mid \filt)\) this null set \emph{must} therefore
    be allowed to depend on \(x\) and \(A\).
}
\[
    \Pr(\rf(x) \in A \mid \filt)(\omega) = \kernel_x(\omega; A).
\]
Corollary \ref{cor: sufficient conditions on collections} then establishes
sufficient conditions on the collection \((\kernel_x)_{x\in \domain}\) to
ensure that \textbf{for all} \(\filt\)-measurable \(X\), all \(A\in \borel(\range)\)
and for \(\Pr\)-almost all \(\omega\in \Omega\),\footref{footnote: evaluation of conditional probability}
\[
    \Pr(\rf(X)\in A \mid \filt)(\omega) = \kernel_{X(\omega)}(\omega; A).
\]
This means that random \(X\) may be treated like deterministic \(x\in\domain\) in the ``formula''
\(\kernel_\bullet(\omega; A)\) for \(\Pr(\rf(\bullet) \in A \mid \cF)(\omega)\).

Above it is implicitly assumed that the function
\(\kernel_{X(\cdot)}(\cdot\,; A)\) is measurable and therefore a well-defined
random variable. To guarantee this, we require that \((\omega, x)\mapsto \kernel_x(\omega; A)\) is a measurable
mapping. This is not guaranteed for arbitrary collections of regular conditional
distributions and warrants the following definition. It represents the
first restriction on the collection \((\kernel_x)_{x\in \domain}\).

\begin{definition}[Joint probability kernels]
    \label{def: joint kernel}
    Let \((\Omega, \cA)\), \((\domain, \cX)\) and \((\range, \cY)\) be
    measurable spaces.  The function \(\kernel\colon (\Omega \times \domain)\times
    \cY \to [0, \infty]\) is called a \emph{joint kernel} for the collection \((\kernel_x)_{x\in \domain}\) 
    of probability kernels \(\kernel_x\colon
    \Omega \times \cY \to [0, \infty]\), if
    \begin{enumerate}[noitemsep]
        \item for all \(x\in \domain\)
        \[ 
            \kernel(\omega, x; A) = \kernel_x(\omega; A) \qquad \forall \omega \in \Omega, A \in \cY,
        \]
        \item 
        the mapping \((\omega, x) \mapsto \kernel(\omega, x; A)\) is measurable
        for all \(A \in \cY\), such that the function \(\kernel\) becomes a
        probability kernel.
    \end{enumerate}
    A joint kernel for a collection of regular conditional distributions is
    called a \emph{joint conditional distribution}.
\end{definition}
\begin{remark}[Strong disintegration]
    \label{rem: strong disintegration}
    A joint kernel is simply a kernel on the product space \(\Omega \times
    \domain\) with target \((\range, \cY)\). This does not mean that standard
    disintegration results yield such a kernel. This is because
    a regular conditional distribution requires a measure. On \(\Omega\) this is
    given by \(\Pr\). But the product space \(\Omega \times \domain\)
    is typically not equipped with a measure. If we take a fixed random
    variable \(X\) in \(\domain\), measurable with respect to \(\filt\), then
    standard disintegration results imply the existence of a regular
    conditional probability distribution 
    \[
        \Pr(\rf(X) \in A \mid \filt)(\omega) = \kernel(\omega; A).
    \]
    But the kernel \(\kernel\) does not take \(X\) as input and is only valid
    for this particular random variable \(X\). Similarly, we could equip the
    product space \(\Omega \times \domain\) with a measure \(\Pr' = \Pr \otimes
    \nu\), where \(\nu\) is a measure on \((\domain, \cX)\). Thus
    \(X(\omega, x) \coloneq x\) is a random variable with respect to \(\Pr'\)
    and \(\Pr_X = \nu\).  Standard disintegration results would then imply, that
    for all \(\filt\otimes \cX\)-measurable \(B\)
    we have
    \[
        \Pr'(\{\rf(X) \in A\} \cap B)
        = \int_B \kernel(\omega, x; A) \Pr(d\omega) \nu(dx).
    \]
    Or in other words,
    \[
        \Pr'(\rf(X) \in A \mid \filt, X=x)(\omega) = \kernel(\omega, x; A)
    \]
    for almost all \(\omega\). However, this kernel \(\kernel\) would depend on \(\nu\)
    and would only be valid for \(\nu\)-almost all \(x\in \domain\). Moreover
    we constructed \(X\) independent from \(\filt\), which required us to add \(X=x\)
    to the conditional part.
    In contrast, a ``joint conditional distribution'' is valid
    \textbf{for all} \(x\in \domain\). And the PIC property we will define in
    the following essentially means that \(\kernel\) is a regular conditional
    distribution \textbf{for all} distributions \(\nu=\Pr_X\). One may therefore call these results
    ``strong disintegration'' results.

    These stronger properties are crucial for applications.
    The explicit relation of the joint conditional distribution to collections
    of conditional distributions is crucial to calculate Bayesian optimizers
    (Example \ref{ex: calculating PI} and \ref{ex: calculate EI}).
    The uniform validity over distributions \(\Pr_X\) of \(X\) is important for Examples \ref{ex: maximizer plugged into PI}
    and \ref{ex: maximizer plugged into EI}.
\end{remark}

Proposition \ref{prop: dependent consistent} establishes the existence and uniqueness of a ``plug-in consistent''
(PIC) joint conditional distribution, where a PIC joint kernel
also admits random input. This proposition also generalizes
the setting of Equation \eqref{eq: measurable evaluation I} to an additional
random variable \(Z\), which will be a helpful tool for our main result.

\begin{definition}[Plug-in consistent -- PIC]
    \label{def: pic simple}
    Let \(Z\) be a random variable in the standard Borel space
	\((E,\borel(E))\) and \(\kernel(\omega, x; B)\) a joint conditional distribution for
    \(\Pr((Z,\rf(x)) \in B\mid \filt)\) indexed by \(x\in \domain\), where
    \((\domain, \cX)\) is a measurable space. Then \(\kernel\) is called
    \emph{plug-in consistent} (PIC), if \textbf{for all \(\filt\) measurable
    random variables \(X\)} in \(\domain\), all \(B \in \borel(E)\otimes
    \borel(\range)\) and for \(\Pr\)-almost all \(\omega\in
    \Omega\),\footref{footnote: evaluation of conditional probability}
	\begin{equation}
        \label{eq: measurable evaluation II}
		\Pr\bigl((Z, \rf(X)) \in B \mid \filt\bigr)(\omega)
		= \kernel(\omega, X(\omega); B).
    \end{equation}
\end{definition}

\begin{prop}[Existence and properties of PIC conditional distributions]
    \label{prop: dependent consistent}
	There exists a PIC joint conditional distribution \(\kernel\)
	for \((Z,\rf(x))\) given \(\filt\) as defined in Definition \ref{def: pic simple}.

    This PIC kernel can be chosen such that \(x\mapsto \kernel(\omega, x; \cdot)\) is continuous with
	respect to the weak topology on the space of measures for all \(\omega\in \Omega\).
	Let \(\tilde\kernel\) be another joint conditional distribution for
	\((Z, \rf(x))\) given \(\filt\), which is continuous for almost all \(\omega\) in
    this sense. Then there
	exists a null set \(N\) such that for all \(\omega\in N^\complement\),
	all \(x\in \domain\) and all borel sets \(B\in \borel(E) \otimes \borel(\range)\)
	\[
		\kernel(\omega,x; B) = \tilde\kernel(\omega,x; B).
	\]
	In particular, \(\tilde\kernel\) is also PIC.
    Finally, let \(g\colon E\times \range \to \real\) be a continuous function with
    \(\E\bigl[\sup_{x\in \domain}|g(Z,\rf(x))| \bigm| \filt\bigr] < \infty\) almost surely, then
    \[
        \rg(x) := \int g(z, y)\, \kernel(\cdot, x; dz \otimes dy) \overset{\as}= \E[g(Z,\rf(x)) \mid \filt]
    \]
    is almost surely continuous in \(x\).
\end{prop}

Corollary \ref{cor: sufficient conditions on collections} provides
sufficient conditions to ensure collections of conditional distributions
\((\kernel_x)_{x\in\domain}\) form a joint kernel and are PIC.

\begin{corollary}[Random dependent evaluation]
    \label{cor: sufficient conditions on collections}
    Let \((\kernel_x)_{x\in \domain}\) be a collection of regular conditional
    distributions for \(\rf(x)\) given \(\filt\), where the random function \(\rf\) 
    is continuous with locally compact, separable
    domain \(\domain\) (e.g.\ \(\real^\dims\)) and Polish codomain \(\range\)
    (e.g.\ \(\real^n\)). Then \textbf{for all} \(\filt\)-measurable \(X\)
    and all \(A\in \borel(\range)\) we have
    \begin{equation}
        \label{eq: measurable evaluation collections}    
        \Pr(\rf(X) \in A\mid \filt)(\omega)
        = \kernel_{X(\omega)}(\omega; A)
    \end{equation}
    for almost all \(\omega\), if
    \begin{enumerate}[label={\normalfont(\roman*)}, noitemsep]
        \item\label{it: joint kernel}
        the map \((\omega, x) \mapsto \kernel_x(\omega; A)\) is
        \textbf{measurable} for all \(A\in \borel(\range)\),
        
        \item\label{it: continuous kernel}
        \(x\mapsto \kernel_x(\omega; \cdot)\)
        is \textbf{continuous} in the weak topology for almost all \(\omega\).
        \end{enumerate}
\end{corollary}
\begin{proof}
    \ref{it: joint kernel} ensures that \(\tilde{\kernel}(\omega, x; A) := \kernel_x(\omega; A)\)
    is a joint kernel and \ref{it: continuous kernel} ensures that \(\tilde{\kernel}\)
    is PIC as it satisfies the requirements in Proposition \ref{prop: dependent consistent}.
\end{proof}

\begin{remark}[Tightness of the result]
    Counterexample \ref{ex: inconsistent joint conditional distribution} shows
    that some continuity in \(x\) is necessary.
\end{remark}

\subsection{General case: Conditionally independent evaluation locations}
\label{sec: random eval in conditional} 

In the introduction we motivated previsible sequences \((X_n)_{n\in \nat}\) with
respect to the filtration \(\filt_n = \sigma(\rf(X_0), \dots, \rf(X_n))\). To
get to our main result we generalize this previsible setting to conditionally
independent evaluation points, admit noisy function evaluation
and also admit starting information.

\paragraph*{Conditional independence}
Sometimes \(X_{n+1}\) is not previsible itself, but sampled from a previsible
distribution.  That is, a distribution constructed from previously seen
evaluations (e.g.\@ Thompson sampling \citep{thompsonLikelihoodThatOne1933}). In
this case,
\(X_{n+1}\) is not previsible, but independent from \(\rf\) conditional on
\(\filt_n\). As introduced in \citet[p.\@
109]{kallenbergFoundationsModernProbability2002}
we denote this by \(X_{n+1} \indep_{\filt_n} \rf\).
In our setting conditional independence is equivalent to\footnote{
    The equivalence only requires \(X_{n+1}\) to be a random variable in a standard Borel space.
} \(X_{n+1} =
h(\xi, U)\) for a measurable function \(h\), a random (previsible) element
\(\xi\) that generates \(\filt_n\) and a standard uniform random variable \(U\)
independent from \((\rf, \filt_n)\) \citep[Prop.\@
6.13]{kallenbergFoundationsModernProbability2002}.

\paragraph*{Noisy evaluations} In many optimization applications only
noisy evaluations of the random objective function \(\rf\) at \(x\) may be
obtained. We associate the noise \(\noise_n\) to the \(n\)-th evaluation
\(x_n\), such that the function \(\rf_n = \rf + \noise_n\) returns the \(n\)-th
observation \(Y_n = \rf_n(x_n)\). While the noise may simply be independent, 
identically distributed constants, observe that this framework allows for much
more general location-dependent noise. The only requirement is that the random
noise functions \(\noise_n\) are continuous, such that the \(\rf_n\) are continuous
functions.

\begin{definition}[Conditionally independent evolution]
    \label{def: cond. indep. evolution}
	The \emph{general conditional independence setting} is given by
    \begin{itemize}[noitemsep]
		\item an underlying probability space \((\Omega, \cA, \Pr)\),
        \item a sub-\(\sigma\)-algebra \(\filt\) (the `initial
        information') with \(W\) a random element such that \(\filt =
        \sigma(W)\),\footnote{
            There always exists such a random element \(W\) since the identity
            map from the measurable space \((\Omega, \cA)\) into \((\Omega,
            \filt)\) is measurable and clearly generates \(\filt\).
        }
        \item A sequence \((\rf_n)_{n\in \nat_0}\) of continuous random
        functions with \(\rf_n\) in \(C(\domain_n, \range_n)\), where the
        domains \((\domain_n)_{n\in \nat_0}\) are locally compact, separable
		metrizable spaces and the codomains \((\range_n)_{n\in \nat_0}\) Polish spaces.
		\item a random variable \(Z\) in a standard Borel space \((E,\borel(E))\)
		(representing an additional quantity of interest).
    \end{itemize}
	A sequence \(X=(X_n)_{n\in \nat_0}\) of random evaluation locations with \(X_n \in \domain_n\) 
	is called a \emph{conditionally independent evolution}, if \(X_{n+1} \indep_{\filt_n} (Z, (\rf_n)_{n\in \nat_0})\)
	for the filtration
	\begin{equation}
        \label{eq: filtration}
		\filt_n \coloneq \sigma(\filt, \rf_0(X_0),\dots, \rf_n(X_n), X_{[0\inter n]})
		\qquad \text{for \(n\ge 0\),} \qquad \filt_{-1} \coloneq \filt,
    \end{equation}
    where \(x_I = (x_i)_{i\in I}\) and we introduce compact notation for \emph{discrete intervals}:
    \[
        \tag{discrete intervals}
        [i\inter j] \coloneq [i,j] \cap \integer,
        \qquad
        [i\inter j) \coloneq [i,j)\cap \integer,
        \qquad
        \text{etc.}
    \]
\end{definition}

In the following we generalize the definition of PIC to this conditional
independence setting before we state our main result.

\begin{definition}[Plug-in consistent -- PIC]
    \label{def: pic}
    Let \(\domain = \prod_{k=0}^n \domain_k\) and let the collection of kernels
    \((\kernel_{x_{[0\inter n]}})_{x_{[0\inter n]} \in \domain}\)
    with
    \[
        \Pr\bigl((Z, \rf_n(x_n)) \in A \mid \filt, (\rf_k(x_k))_{k\in [0\inter n)}\bigr)
        \overset\as= \kernel\bigl(W, (\rf_k(x_k))_{k\in[0\inter n)}, x_{[0\inter n]}; A\bigr)
    \]
    form a joint kernel \(\kernel\) (Definition \ref{def: joint kernel}). Then this joint
    kernel \(\kernel\) is called \emph{plug-in consistent} (PIC), if \textbf{for all} conditionally independent
    evolutions \((X_k)_{k\in \domain}\) and their filtration \(\filt_n\) as defined in \eqref{eq: filtration},
    we have
    \[
        \Pr\bigl((Z, \rf_n(X_n)) \in A \mid \filt_{n-1}, X_n\bigr)
        \overset\as=
        \kernel\bigl(W, (\rf_k(X_k))_{k\in [0\inter n)}, X_{[0\inter n]}; A\bigr).
    \]
\end{definition}
\begin{remark}[This is a generalization]
    Note that PIC as defined above is a generalization of
    Definition \ref{def: pic simple} with \(n=0\) and \(X=X_0\).
\end{remark}

In our main result we allow for random evaluations in the conditional.
This is significantly harder to prove than random evaluation locations
in the dependent part of the conditional distribution.
The main ingredient of the proof is a general result that allows
moving plug-in consistency (PIC) from the dependent to the conditional part
(Proposition \ref{prop: consistency shuffle}). We iteratively combine
this result with Proposition \ref{prop: dependent consistent} that only allows
random variables in the dependent part.

\begin{theorem}[Conditionally independent sampling]
	\label{thm: conditionally independent sampling}
	Assume the general conditional independence setting (Definition \ref{def: cond. indep. evolution}).
    \begin{enumerate}[label={(\roman*)}]
        \item\label{it: consistent without dependent variable}
        \textbf{Without the dependent variable, every joint conditional distribution is PIC.}
        That is, let \(\kernel\) be a joint conditional distribution for \(Z\) given
        \((\filt, \rf_0(x_0), \dots, \rf_n(x_n))\), that is for all \(x_{[0\inter n]}\in
        \prod_{k=0}^n\domain_k\) and \(A\in\borel(E)\)
        \[
            \Pr\bigl(Z\in A \mid \filt, \rf_0(x_0), \dots, \rf_n(x_n)\bigr)
            \overset{\as}= \kernel\bigl(W, \rf_0(x_0), \dots, \rf_n(x_n), x_{[0\inter n]}; A\bigr).
        \]
        Then for all conditionally independent evolutions
        \((X_k)_{k\in \nat_0}\) and \(A\in\borel(E)\)
        \[
            \Pr\bigl(Z\in A \mid \filt_n\bigr)
            \overset\as=
            \kernel\bigl(W, \rf_0(X_0), \dots, \rf_n(X_n), X_{[0\inter n]}; A\bigr).
        \]
        \item\label{it: consistent with dependent variable}
        \textbf{With the dependent variable, continuity is sufficient for joint conditional
        distributions to be PIC.} That is,
		let the kernel \(\kernel\) be a joint conditional distribution for
        \((Z,\rf_n(x_n))\) given \((\filt, (\rf_k(x_k))_{k\in[0\inter n)})\)
        such that
        \[
            x_n \mapsto \kernel(y_{[0\inter n)}, x_{[0\inter n]}; \,\cdot\,)
        \]
        is continuous in the weak topology for all \(x_{[0\inter
        n)}\in \domain^n\) and all \(y_{[0\inter n)}\in \range^n\).
		Then for all conditionally independent evolutions \((X_k)_{k\in \nat_0}\)
		and measurable sets \(B\in\borel(E)\otimes \borel(\range)\)
        \[
            \Pr\bigl((Z,\rf_n(X_n)) \in B \mid \filt_{n-1}, X_n\bigr)
            \overset\as=
            \kernel\bigl(W, (\rf_k(X_k))_{k\in [0\inter n)}, X_{[0\inter n]}; B\bigr).
        \]
    \end{enumerate}
\end{theorem}

\begin{remark}[Comparison with dependent evaluation]
    We highlight that continuity of the kernel is only required for the case
    where function evaluations are dependent variables.
    However, while Proposition \ref{prop: dependent consistent} shows the
    existence of a joint kernel that admits random evaluation, the existence of such a joint conditional
    distribution is not proven here. In the Gaussian case this is not a problem,
    since the joint conditional distribution is known explicitly (see Section~\ref{sec: gaussian conditionals}).
\end{remark}

\begin{corollary}[Gaussian case]
    \label{cor: gaussian case}
    In the general conditional independence setting (Definition \ref{def: cond. indep. evolution})
    further assume that the sequence \((\rf_n)_{n\in \nat_0}\) consists of joint
    Gaussian random functions. Then the canonical Gaussian conditional
    distribution (Definition \ref{def: canonical Gaussian dist}) for
    \(\rf(x,n) = \rf_n(x)\) is PIC.
\end{corollary}
\begin{proof}

    The canonical Gaussian conditional distribution is a joint conditional distribution (Remarks \ref{rem: Gaussian join kernel})
    that also satisfies the continuity requirement in Theorem \ref{thm:
    conditionally independent sampling} \ref{it: consistent with dependent
    variable} (Remark \ref{rem: continuity of gaussian conditional}).
\end{proof}

\subsection{Counterexample}
\label{sec: counterexample}

While it is unlikely that anyone would question the measurability requirement in
Corollary \ref{cor: sufficient conditions on collections}, the
continuity requirement may seem strange. Especially, as we do not appear to
require it for the evaluation of functions in the
conditional (see Theorem \ref{thm: conditionally independent sampling}).
The following counterexample illustrates that such a continuity
requirement is indeed necessary for the PIC property.

\begin{counterexample}[Joint kernel that is not PIC] 
    \label{ex: inconsistent joint conditional distribution}
    In this example we show that \(X\) measurable with respect to \(\filt\) is
    not always sufficient to treat \(X\) as deterministic in \(\E[\rf(X) \mid
    \filt]\).

    For this purpose consider a standard normal random variable
    \(Y\sim \normal(0,1)\) independent of a standard uniform random variable
    \(U\sim \uniform(0,1)\). We define an almost constant Gaussian random function
    \(\rf(x) := Y\) and define \(\filt := \sigma(U)\). Since the
    conditional expectation is only well-defined up to null sets we have
    \[
        \E[\rf(x) \mid U] \overset{\as}= \E[Y] \ind_{U \neq x} =: g(x).
    \]
    However we have
    \[
        g(U) = 0 \neq \E[Y] \overset{\as}= \E[\rf(U) \mid U].
    \]
    While the formula \(g\) is therefore a valid formula for \(\E[\rf(x)\mid
    U]\) for all deterministic \(x\), it is not valid for random \(X=U\). Even
    though \(X\) is measurable with respect to \(\filt=\sigma(U)\).
    Similarly one could construct \(g\) to not be a measurable
    function. In this case \(g(U)\) would not even be a valid random variable.
    This further justifies \ref{it: joint kernel} in Corollary \ref{cor:
    sufficient conditions on collections}.
\end{counterexample}

\begin{remark}[Connection to regular conditional probability]
    In the Counterexample above we created different adversarial null sets for
    each \(x\), which joined together in the case of the random \(X=U\) to
    break the formula. Similar counterexamples motivate the construction of
    the \emph{regular conditional probability kernel}\footnote{
        The conditional
        probability \(\Pr(X\in A \mid \filt) = \E[\ind_A \mid \filt]\) for a random
        variable \(X\), \(\sigma\)-algebra \(\filt\) and measurable sets \(A\)
        is only well defined up to a null set. It is therefore ex ante impossible to
        ensure that \(A \mapsto \Pr(X \in A \mid \filt)\) is a well-defined measure.
        This is because the null-set may depend on \(A\) and, if adversarially
        selected, their union may no longer be a null set.
    } \citep[e.g.][Def.\@ 8.28]{klenkeProbabilityTheoryComprehensive2014}.
    Just like regular conditional probabilities result in a narrower definition of
    a ``sensible'' conditional probability, PIC joint conditional distributions
    further narrow the set of conditional distributions down.
\end{remark} 	\section{Applications}
\label{sec: applications}

In the case of Gaussian random functions \(\rf\) for example,
\((\rf(x_0), \dots, \rf(x_n))\) is a multivariate Gaussian vector
with well known posterior distribution \(\rf(x_n)\) given \((\rf(x_0),\dots,
\rf(x_{n-1}))\) when the evaluation locations are deterministic. But \(\rf(X)\)
is not necessarily Gaussian if \(X\) is random\footnote{
    consider \(X= \argmin_{x\in K}\rf(x)\) for some compact set \(K\subseteq
    \domain\).
} and the calculation of conditional distributions hence becomes much more
difficult. Treating previsible inputs as deterministic ensures the
calculation is feasible but it requires a theoretical foundation.

\subsection{Maximal probability of improvement}

To illustrate the necessity of our results, we consider a simple optimization
procedure from Bayesian optimization known as ``maximal probability of
improvement'' (PI) \citep[e.g.][Sec.\ 7.5]{garnettBayesianOptimization2023}.

Let \((\rf(x))_{x\in \domain}\) be a Gaussian random function with values in
\(\real\). The goal is to find a maximum of \(\rf\). To this end, new
evaluation locations are chosen according to the rule
\begin{equation}
    \label{eq: PI}
    \tag{PI}
    X_n := \argmax_{x_n\in\domain}\Pr\Bigl(\rf(x_n) > \underbrace{\max_{i=0,\dots, n-1} \rf(X_i) + \epsilon}_{\eqcolon\eta} \mid \filt_{n-1}\Bigr),
\end{equation}
where \(\epsilon>0\) is a minimum improvement and \(\filt_n = \sigma(\rf(X_0), \dots, \rf(X_n))\).
 And we assume a deterministic starting location \(X_0= x_0\). We will address
the measure-theoretic subtleties 
of maximization in Subsection \ref{sec: subtleties of maximization}, let us first consider whether or not we can even
obtain an explicit expression for the probability of improvement.

\subsubsection{Calculating the probability of improvement}

\paragraph*{Deterministic case}
Computing \ref{eq: PI} explicitly would be straightforward if the \(X_0,\dots, X_{n-1}\) were
deterministic \(x_0, \dots, x_{n-1}\). In this case \((\rf(x_0), \dots, \rf(x_n))\) is
a multivariate Gaussian random vector. Consequently, the conditional distribution of
\(\rf(x_n)\) given \((\rf(x_0), \dots, \rf(x_{n-1}))\) is again Gaussian
(see Lemma~\ref{lem: gaussian conditional}).
Let \(\Phi\) denote the cumulative distribution function of the standard
normal distribution \(\normal(0,1)\).  Since \(\eta\) is \(\filt_{n-1}\) measurable we have
\[
    \Pr\bigl(\rf(x_n) > \eta \mid \filt_{n-1}\bigr)
    = 1- \Phi(\tfrac{\eta - \mu_n}{\sigma_n}),
\]
where the posterior mean and variance,
\begin{align*}
    \mu_n
    = \mu_n(x_0,\dots, x_n, \rf(x_0), \dots, \rf(x_{n-1}))
    \quad \text{and}\quad
    \sigma_n^2
    = \sigma_n^2(x_0,\dots, x_n),
\end{align*}
are given explicitly in Lemma \ref{lem: gaussian conditional}.
Since \(\Phi\) is monotonically increasing this results 
in the optimization problem
\begin{align}
    X_n &= \argmax_{x_n}
    \Pr\bigl(\rf(x_n) > \eta \mid \filt_{n-1}\bigr)
    = \argmax_{x_n} \frac{\mu_n - \eta}{\sigma_n},
    \label{eq: PI objective}
\end{align}
which can be numerically maximized using the explicit formulas
in Lemma \ref{lem: gaussian conditional}.

\paragraph*{Heuristic in the random case}
With the heuristic motivation that the \(X_k\) are `deterministic' conditioned
on \(\filt_{k-1}\) (``previsible''), the same procedure is used when the
\(X_0, \dots, X_{n-1}\) are \emph{not} deterministic but selected by this
process. The correctness of this procedure is often treated as self-evident
\citep[e.g.][Lemma 5.1, p.\
3258]{srinivasInformationTheoreticRegretBounds2012}.
A proof of this conjecture is non-trivial but a direct result of Corollary \ref{cor: gaussian
case}.

\begin{example}[Calculating the probability of improvement]
    \label{ex: calculating PI}
    Let \(\rf\) be a continuous Gaussian random function, \(X_k\) measurable
    with respect to \(\filt_{k-1} = \sigma(\rf(X_0), \dots, \rf(X_{k-1}))\)
    and \(X_0 = x_0\) deterministic.
    Then, if the probability of improvement
    \[ 
        \Pr\Bigl(\rf(x_n) > \max_{k\in[0\inter n)} \rf(X_k) + \epsilon \mid \filt_{n-1}\Bigr)
    \]
    is calculated using the canonical Gaussian conditional distribution (Definition \ref{def:
    canonical Gaussian dist})
    the random evaluation locations \(X_k\) may be treated like deterministic locations \(x_k\).
\end{example}
\begin{proof}
    Define the continuous function
    \[
        h(y_{[0\inter n]}) \coloneq y_n -  \max_{k\in[0\inter n)} y_k.
    \]
    Using \(\rf_n(x_{[0\inter n]}) \coloneq (\rf(x_0), \dots, \rf(x_n))\) the desired probability is then clearly 
    \[
        \Pr\Bigl(
            \rf_n(X_{[0\inter n)}, x_n) \in h^{-1}((\epsilon, \infty))
            \mid \rf(X_0), \dots, \rf(X_{n-1})
        \Bigr).
    \]
    Since the \((\rf_k)_{k\le n}\) with \(\rf_k \coloneq \rf\) for \(k< n\) are
    joint Gaussian, the canonical Gaussian conditional distribution (Definition
    \ref{def: canonical Gaussian dist}) for deterministic locations forms a
    joint kernel \(\kernel\) (Remark \ref{rem: Gaussian join kernel}) with 
    \[
        \Pr\Bigl(\rf_n(\tilde{x}_n) \in A \mid \rf(x_0), \dots, \rf(x_{n-1})\Bigr)
        = \kernel(\rf(x_0), \dots, \rf(x_{n-1}), x_{[0\inter n)}, \tilde x_n; A),
    \]
    We may then apply Corollary \ref{cor: gaussian case} to obtain that
    we have for \(\tilde X_n = (X_{[0\inter n)}, x_n)\)
    \[
        \Pr\Bigl(\rf(x_n) > \max_{k\in[0\inter n)} \rf(X_k) + \epsilon \mid \filt_{n-1}\Bigr)
        = \kernel(\rf(X_0), \dots, \rf(X_{n-1}), X_{[0\inter n)}, \tilde X_n; A).
    \]
    Or in other words, we may treat the inputs \(X_0, \dots, X_{n-1}\) like deterministic parameters
    in the calculation.
\end{proof}

\subsubsection{Subtleties of maximization}
\label{sec: subtleties of maximization}

Assuming there exists a continuous version of \(x \mapsto \Pr(\rf(x) > \eta \mid \filt_{n-1})\)
and \(\domain\) is compact, the maximizer \(X_n\) in \ref{eq: PI} exists.
And it is perhaps reasonable to expect
\begin{equation}
    \label{eq: maximizer plugged in is max}
    \Pr\bigl(\rf(X_n) > \eta \mid \filt_{n-1}\bigr)
    \overset\as= \max_{x_n \in \domain} \Pr\bigl(\rf(x_n) > \eta \mid \filt_{n-1}\bigr),
\end{equation}
since \(X_n\) is defined as the maximizer. However, the supremum of a random
function is well-known to be a subtle measure-theoretic object that is only well
defined if the random function is separable. The conditional probabilities are
such random functions in \(x_n\) (assuming they are formed using a joint
kernel for measurability).  The following illustrative example explains why we
cannot expect \eqref{eq: maximizer plugged in is max} in general.
\begin{example}[Maximization of conditional expectations]
    Let \(U\sim \uniform([0,1])\) and \(Y \sim \normal(0,1)\) be independent
    random variables and define \(\rf(x) = Y\) to be the constant Gaussian
    random function over \(x\in [0,1]\). Since the conditional expectation is
    only well defined up to null sets we clearly have
    for \(\gamma \gg \E[Y]\)
    \[
        \E[\rf(x) \mid U] \overset\as= \E[Y] + \ind_{x=U} \gamma =: g(x),
    \]
    as \(\{x=U\}\) is a null set. Therefore we have
    \[
        \max_{x\in [0,1]} g(x) = \gamma, \qquad \argmax_{x\in [0,1]} g(x) = U.
    \]
    However this does not translate back
    \[
        \E[\rf(U) \mid U] \overset{\as}= \E[Y] \neq g(U) = \gamma.
    \]
    Since \(g(x)\) is a valid version of \(\E[\rf(x) \mid U]\) the term
    \[
        \sup_{x\in [0,1]} \E[\rf(x) \mid U] 
    \]
    is not well defined if we admit such \(g\). However, if we require
    \(x\mapsto \E[\rf(x) \mid U]\) to be separable (e.g.\ continuous), then
    the supremum is well defined.
\end{example}

We now may want to prove \eqref{eq: maximizer plugged in is max} for separable versions of
the conditional distribution with Proposition \ref{prop: dependent consistent}.
For this we require the following about the PIC joint kernel for all \(x_n \in \domain\)
and all \(\eta \in \real\)
\begin{equation}
    \label{eq: no atom}    
    \Pr\bigl(\rf(x_n) = \eta \mid \filt_{n-1}\bigr) \overset\as= \kernel(\cdot, x_n; \{\eta\}) = 0.
\end{equation}
This is generally the case for Gaussian random functions, since
the value
\[
    \eta = \max_{i=0,\dots, n-1} \rf(X_i) + \epsilon
\]
is not achieved at any of the locations \(X_i\) for \(\epsilon>0\) and all other
points typically have positive remaining variance.
\begin{example}[Maximizer plugged into PI]
    \label{ex: maximizer plugged into PI}
    Assuming \eqref{eq: no atom} we have
    \begin{align*}
        \Pr\bigl(\rf(X_n) > \eta \mid \filt_{n-1}\bigr)
        \overset\as&= \max_{\substack{X \text{ \(\filt_{n-1}\)-meas.}\\\text{r.v. in \(\domain\)}}} \Pr\bigl(\rf(X) > \eta \mid \filt_{n-1}\bigr).
        \\
        \overset\as&= \max_{x_n \in \domain} \Pr\bigl(\rf(x_n) > \eta \mid \filt_{n-1}\bigr).
    \end{align*}
\end{example}
\begin{proof}
    Using \(h(y_{[0\inter n]}) = y_n - \max_{k\in[0\inter n)} y_k\) and
    \(\rf_n(x_{[0\inter n]}) := (\rf(x_0), \dots, \rf(x_n))\) again, we have
    by Proposition \ref{prop: dependent consistent} that there exists a kernel
    with
    \begin{align*}
        \Pr(\rf_n(X_{[0\inter n) }, x_n) \in A \mid \filt_{n-1})(\omega)
        &= \kernel(\omega, X_{[0\inter n)}, x_n; A),
    \end{align*}
    where \(\tilde x \mapsto \kernel(\omega, \tilde x; \cdot)\) is continuous
    in the weak topology of measures. Since \((\epsilon, \infty)\) is a continuity set, i.e.\ \(\Pr(\rf_n(X_{[0\inter n) }, x_n) = \epsilon \mid \filt_{n-1}) = 0\)
    by \eqref{eq: no atom}, we have by the Portmanteau theorem \citep[e.g.][Thm.\@ 13.16]{klenkeProbabilityTheoryComprehensive2014}
    that the following function is continuous
    \[
        x_n \mapsto \kernel(\omega, X_{[0\inter n)}, x_n; (\epsilon, \infty)).
    \]
    This function is thereby a separable version of PI and thus
    \begin{align*}
        \max_{x_n \in \domain} \Pr\bigl(\rf(x_n) > \eta \mid \filt_{n-1}\bigr)
        \overset\as&= \max_{x_n \in \domain} \kernel(\cdot, X_{[0\inter n)}, x_n; [\epsilon, \infty))
        \\
        &\ge \kernel(\cdot, X_{[0\inter n)}, X_n; [\epsilon, \infty)) 
        \\
        \overset\as&= \Pr\bigl(\rf(X_n) > \eta \mid \filt_{n-1}\bigr)
    \end{align*}
    for any \(X_n\) that is \(\filt_{n-1}\)-measurable and takes values in \(\domain\).
    Taking the maximum over such \(X_n\) results in the second equality since
    deterministic \(x_n\) are special cases. Since
    \[
        X_n = \argmax_{x_n}
        \Pr\bigl(\rf(x_n) > \eta \mid \filt_{n-1}\bigr)
        \overset\as= \argmax_{x_n} \kernel(\cdot, X_{[0\inter n)}, x_n; [\epsilon, \infty))
        \eqcolon \tilde X_n,
    \]
    we clearly also have the first equality
    \begin{align*}
        \Pr\bigl(\rf(X_n) > \eta \mid \filt_{n-1}\bigr)
        \overset\as&= \kernel(\cdot, X_{[0\inter n)}, X_n; [\epsilon, \infty))
        \\
        \overset\as&= \kernel(\cdot, X_{[0\inter n)}, \tilde X_n; [\epsilon, \infty))
        \\
        &= \max_{x_n \in \domain} \kernel(\cdot, X_{[0\inter n)}, x_n; [\epsilon, \infty))
        \\
        \overset\as&= \max_{x_n \in \domain} \Pr\bigl(\rf(x_n) > \eta \mid \filt_{n-1}\bigr).
        \qedhere
    \end{align*}
\end{proof}

\begin{remark}
    It may be possible to get rid of the assumption \eqref{eq: no atom} if
    one can combine the continuity of the kernel in the weak sense with the fact
    that the function \(x \mapsto \Pr(\rf(x) > \eta \mid \filt_{n-1})\)  for the fixed set \([\epsilon, \infty)\) must be
    separable.
\end{remark}

\subsection{Maximal expected improvement}

Another popular method in Bayesian optimization is the maximization
of the \emph{expected improvement} (EI) \citep[e.g.][Sec.\ 7.3]{garnettBayesianOptimization2023}.
Here the next evaluation location is chosen as
\[
    X_n = \argmax_{x_n\in \domain} \E\Bigl[ \bigl(\rf(x_n) - \rf_n^*\bigr)_+ \mid \filt_{n-1}\Bigr].
\]
for \(\rf_n^* = \max_{k\in[0\inter n)} \rf(X_k)\)
with \(x_+ = \max\{0, x\}\), assuming \(\sup_{x\in \domain} \E|\rf(x)| < \infty\).
\begin{example}[Expected improvement with previsible inputs]
    \label{ex: calculate EI}
    Let \(\rf\) be a continuous Gaussian random function and \(X_k\)
    measurable with respect to the filtration \(\filt_{k-1} = \sigma(\rf(X_0), \dots, \rf(X_{k-1}))\)
    and \(X_0 = x_0\) deterministic.
    Then, the expected improvement
    \[
        \E\Bigl[ \bigl(\rf(x_n) - \rf_n^*\bigr)_+ \mid \filt_{n-1}\Bigr]
    \]
    may be calculated using the canonical conditional distribution (Definition \ref{def: canonical Gaussian dist})
    by treating the random evaluation locations \(X_k\) like deterministic
    locations \(x_k\).
\end{example}
\begin{proof}
    Using the continuous function
    \[
        h(y_{[0\inter n]}) = (y_n - \max_{k\in[0\inter n)} y_k)_+
    \]
    and \(\rf_n(x) = (\rf(x_0), \dots, \rf(x_n))\) the expected improvement
    can be written as
    \[
        \E\Bigl[ h \circ \rf_n(X_{[0\inter n)}, x_n) \mid \filt_{n-1}\Bigr]
        = \int h(y) \, \Pr(\rf_n(X_{[0\inter n)}, x_n) \in dy\mid \filt_{n-1})
    \]
    Application of Corollary \ref{cor: gaussian case} 
    to \(\Pr(\rf_n(X_{[0\inter n)}, x_n) \in dy\mid \filt_{n-1})\) analogous
    to the proof of Example \ref{ex: calculating PI} yields the claim.
\end{proof}
Without the assumption \(\E \sup_{x\in \domain} |\rf(x)| < \infty\) it is
already difficult to justify the existence of a continuous conditional
expectation \(x \mapsto \E[\rf(x) \mid \filt]\). The following
assumption in the example on maximization is therefore very natural.
\begin{example}[Maximizers plugged into expected improvement]
    \label{ex: maximizer plugged into EI}
    Assuming \(\E[\sup_{x\in \domain} |\rf(x)|] < \infty\)
    we have 
    \begin{align*}
        \E\Bigl[ \bigl(\rf(X_n) - \rf_n^*\bigr)_+ \mid \filt_{n-1}\Bigr]
        &= \max_{\substack{X \text{ \(\filt_{n-1}\)-meas.}\\\text{r.v. in \(\domain\)}}} \E\Bigl[ \bigl(\rf(X) - \rf_n^* \bigr)_+ \mid \filt_{n-1}\Bigr]
        \\
        &= \max_{x_n \in \domain} \E\Bigl[ \bigl(\rf(x_n) - \rf_n^*\bigr)_+ \mid \filt_{n-1}\Bigr].
    \end{align*}
\end{example}
\begin{proof}
    Using \(h(y_{[0\inter n]}) = (y_n - \max_{k\in[0\inter n)} y_k)_+\) and \(\rf_n\) again, we have
    \[
        \E[\sup_{x\in \domain} |h\circ \rf_n(x_{[0\inter n]}) | \mid \filt_{n-1}]
        \le (n+1)\E\Bigl[\sup_{x\in \domain}|\rf(x)| \mid \filt_n \Bigr] < \infty
    \]
    almost surely. For the kernel \(\kernel\) from Proposition \ref{prop:
    dependent consistent} with
    \[
        \Pr(\rf_n(x_{[0\inter n) }, x_n) \in A \mid \filt_{n-1}) = \kernel(\cdot, x_{[0\inter n)}, x_n; A)
    \]
    we thereby know that
    \[
        H(x_n) \coloneq \int h(y) \, \kernel(\cdot, X_{[0\inter n)}, x_n; dy)
        \overset\as= \E\Bigl[ \bigl(\rf(x_n) - \rf_n^* \bigr)_+ \mid \filt_{n-1}\Bigr]
    \]
    is an almost surely continuous version of the expected improvement into
    which we may plug random \(X_n\) that are measurable with respect to
    \(\filt_{n-1}\).  The rest of the proof is then analogous to the proof of
    Example \ref{ex: maximizer plugged into PI}.
\end{proof}

	\section{Proofs}
\label{sec: proofs}

\subsection{Simplified case: Dependent evaluation}

\begin{proof}[Proof of Proposition \ref{prop: dependent consistent}]
	Observe that \(E\times C(\domain, \range)\) is a standard borel space since \(C(\domain,
	\range)\) is Polish (Theorem~\ref{thm: topology of continuous functions}). There
	therefore exists a regular conditional probability distribution
	\(\kernel_{Z,\rf \mid \filt}\) \citep[e.g.][Thm.\@ 6.3]{kallenbergFoundationsModernProbability2002}.
	Using this probability kernel, we define the kernel
	\[
		\kernel(\omega, x; B)
		\coloneq \int \ind_B(z, e(f,x)) \kernel_{Z, \rf\mid \filt}(\omega; dz \otimes df)
	\]
	which is a measure in \(B\in \borel(E)\otimes \borel(\range)\) by linearity of the integral and
	monotone convergence. We therefore only need to prove it is measurable in
	\((\omega, x)\in \Omega\times\domain\) to prove it is a probability kernel.
	This follows from measurability of the evaluation function \(e\)
	(Theorem~\ref{thm: topology of continuous functions}) and the application of Lemma 14.20
	by \citet{klenkeProbabilityTheoryComprehensive2014}
	to the probability kernel \(\tilde\kernel(\omega, x; A) \coloneq \kernel_{Z,\rf\mid \filt}(\omega; A)\)
	in the equation above. By `disintegration'
	\citep[e.g.][Thm~6.4]{kallenbergFoundationsModernProbability2002}
	this probability kernel is moreover a regular conditional version of
	\(\Pr((Z,\rf(X))\in B\mid \filt)\) for all \(\filt\)-measurable \(X\), i.e. for all \(B\in \borel(E)\otimes
	\borel(\range)\) and for \(\Pr\)-almost all \(\omega\)
	\begin{align*}
		\Pr\bigl((Z, \rf(X)) \in B \mid \filt\bigr)(\omega)
		\overset{\text{disint.}}&= \int \ind_B(z, e(f,X(\omega))) \kernel_{Z,\rf\mid \filt}(\omega; dz \otimes df)
		\\
		\overset{\text{def.}}&= \kernel(\omega, X(\omega); B).
	\end{align*}
	The kernel thereby satisfies \eqref{eq: measurable evaluation II}.
	
	For continuity observe that we have \(\lim_{x\to y} (z,f(x)) = (z,f(y))\)
	for any \(f\in C(\domain, \range)\). For open \(U\) this implies
	\[
		\liminf_{x\to y} \ind_U(z,f(x)) \ge \ind_U(z,f(y)),
	\]
	because if \((z,f(y)) \in U\), then eventually \((z,f(x))\) in \(U\) due to
	openness of \(U\). An application of Fatou's lemma \citep[e.g.][Thm.\@ 4.21]{klenkeProbabilityTheoryComprehensive2014}
	yields for all open \(U\)
	\[
		\liminf_{x\to y}\kernel(\omega, x; U)
		\ge \int \liminf_{x\to y} \ind_U(z,f(x)) \kernel_{Z,\rf\mid \filt}(\omega; dz\otimes df)
		\ge \kernel(\omega, y; U).
	\]
	And we can conclude weak convergence by the Portmanteau theorem
	\citep[Thm.\@ 13.16]{klenkeProbabilityTheoryComprehensive2014} since \(E\times \range\)
	is metrizable.

	Let \(\tilde\kernel\) be another joint probability kernel that is continuous almost
	surely. We will assume it is continuous for all \(\omega\) without loss of generality
	in the following and assume that the null set of discontinuity is tacitly joined
	with the null set we construct.
	Since \(E\times \range\) is second countable, there is a countable base
	\(\{U_n\}_{n\in \nat}\) of its topology, which generates the Borel
	\(\sigma\)-algebra \(\borel(E) \otimes \borel(\range)\). And since \(\domain\)
	is separable, it has a countable dense subset \(Q\). 
	There must therefore exist a null set \(N\) such that
	\[
		\kernel(\omega, q; U_n)
		= \tilde\kernel(\omega, q; U_n),
		\qquad \forall \omega \in N^\complement,\ n\in \nat,\ q\in Q,
	\]
	because both kernels are regular conditional version of \(\Pr(Z,\rf(q)\in U_n; \cG)\)
	and the union over \(\nat\times Q\) is a countable union. Since
	\(\{U_n\}_{n\in \nat}\) generates the \(\sigma\)-algebra, we deduce for all
	\(\omega \in N^\complement\) and all \(q\in Q\) that \(\kernel(\omega,q;
	\cdot) = \tilde\kernel(\omega, q; \cdot)\).
	As \(Q\) is dense in \(\domain\) we have by continuity of the joint kernels
	for all \(\omega\in N^\complement\) and all \(x\in\domain\)
	\[
		\kernel(\omega,x; \cdot) 
		= \tilde\kernel(\omega, x; \cdot).
	\]
	Let \(g\colon E\times \range \to \real\) be a continuous function with
	\(\E[\sup_{x}|g(Z,\rf(x))| \mid \filt] < \infty\) almost surely.
	Then for \(\rg(x) = \int g(z, y) \kernel(\cdot, x; dz\otimes dy)\) we have
	\begin{align*}
		|\rg(x) - \rg(x_0)|(\omega)
		&= \Bigl|
			\int g(z, y) \kernel(\omega, x; dz \otimes dy)
			- \int g(z, y) \kernel(\omega, x_0; dz \otimes dy)
		\Bigr|
		\\
		&= \Bigl|
			\int g(z, e(f, x)) - g(z, e(f, x_0)) \kernel_{Z, \rf\mid \filt}(\omega; dz \otimes df)
		\Bigr|
		\\
		&\le 
			\int \bigl|g(z, f(x)) - g(z, f(x_0))\bigr|
			\kernel_{Z, \rf\mid \filt}(\omega; dz \otimes df)
	\end{align*}
	using the definition of \(\kernel\) via \(\kernel_{Z,\rf\mid \filt}\) for the second equation. Up to a null set \(N\)
	we further have for all \(\omega \in N^\complement\)
	\[
		\int \sup_{x} |g(z, f(x))| \kernel_{Z, \rf \mid \filt}(\omega; dz\otimes df)
		= \E\Bigl[\sup_{x}|g(Z, \rf(x))| \Bigm| \filt\Bigr](\omega) < \infty.
	\]
	Since \(|g(z, f(x)) - g(z, f(x_0))| \le 2\sup_x |g(z,\rf(x))|\) we thus have by
	dominated convergence and continuity of \(g\) and \(\rf\) that \(\rg(x) \to
	\rg(x_0)\) for \(x\to x_0\).
\end{proof}

\begin{remark}[Possible generalization]
	\label{rem: existence of consistent joint cond. dist}
	Note that for the existence of a joint conditional distribution that
	satisfies \eqref{eq: measurable evaluation II}
	we only require a regular conditional distribution for \((Z,\rf)\) given
	\(\filt\) to exist and measurability of the evaluation map \(e\). This part
	of the result can therefore be made to hold with greater generality.
\end{remark}

\subsection{General case: Conditionally independent evaluation locations}

To prove our main result (Theorem \ref{thm: conditionally independent sampling})
we will move the PIC property that we can obtain by Proposition \ref{prop:
dependent consistent} for function values in the dependent part to the
conditional part. The general proposition we will prove in the following
is the key tool for this purpose. We change the notation of the random objects
to highlight that we no longer assume continuity or make assumptions about the
domain and range.
This proposition essentially states:

\begin{quote}
	Let \(\xi_1, \xi_3\) be random elements and \(\xi_2 = (\xi_2^y)_{y}\) be a collection of random elements indexed by \(y\).
	If there exists a PIC joint conditional distribution \(\magenta{\kappa_{3,2\mid 1}}\)
	for \((\xi_3, \xi_2^y)\) given \(\xi_1\), then \textbf{any} joint conditional distribution
	\(\teal{\kappa_{3\mid 2,1}}\) for \(\xi_3\) given \((\xi_2^y, \xi_1)\) is
	PIC.
\end{quote}

The ``joint'' in ``joint conditional distribution'' refers to the index \(y\)
of \(\xi_2^y\). The formal statement of this claim follows:

\begin{prop}[Consistency shuffle]\label{prop: consistency shuffle}
	Let \((\Omega, \cA, \Pr)\) be a probability space and let \(\xi_1,
	(\xi_2^{y})_{y\in D},\xi_3\) be random elements in the measurable spaces
	\((E_i, \cE_i)\), \(i\in\{1,2,3\}\), where \(\xi_2\) is indexed by the
	measurable domain \((D, \mathcal D)\).  Assume \(\xi_2^{g(\xi_1)}\) is
	measurable as a random element in \(\omega\in \Omega\) for all measurable functions
	\(g\colon E_1 \to D\),
	and assume there exists a PIC joint conditional distribution
	\(\magenta{\kappa_{3,2\mid 1}}\) for \(\xi_3, \xi_2^y\) given \(\xi_1\). That is 
	for all \(A\in \cE_3\otimes \cE_2\) and all measurable functions \(g\colon E_1 \to D\)
	\begin{equation}
		\label{eq: consistent outer kernel}
		\Pr\bigl(\xi_3, \xi_2^{g(\xi_1)} \in A \mid \xi_1\bigr)
		\overset{\as}= \magenta{\kappa_{3,2\mid 1}}(\xi_1, g(\xi_1); A).
	\end{equation}
	If there exists a joint conditional probability kernel
	\(\teal{\kappa_{3\mid1,2}}\) for \(\xi_3\) given \(\xi_1 ,\xi_2^{y}\) such
	that for all \(y\in D\) and \(A_3\in \cE_3\)
	\begin{equation}
		\label{eq: joint inner kernel}
		\Pr(\xi_3\in A_3 \mid \xi_1, \xi_2^{y})
		\overset{\as}= \teal{\kappa_{3\mid 2, 1}}(\xi_1, \xi_2^{y}, y; A_3),
	\end{equation}
	then \(\teal{\kappa_{3\mid 2, 1}}\) is PIC, that is we have for all \(A_3\in \cE_3\)
	and measurable \(g\colon E_1 \to D\)
	\[
		\Pr(\xi_3 \in A_3 \mid \xi_1, \xi_2^{g(\xi_1)})
		\overset{\as}= \teal{\kappa_{3\mid2,1}}\bigl(\xi_1, \xi_2^{g(\xi_1)}, g(\xi_1); A_3\bigr).
	\]
\end{prop}
\begin{remark}[Possible generalization]
	Note that we keep \(g\) fixed throughout the proof. So if \(\magenta{\kappa_{3,2\mid 1}}\) is PIC
	only for a specific \(g\), then we also obtain consistency of \(\teal{\kappa_{3\mid2,1}}\)
	only for this specific function \(g\). For consistency of \(\teal{\kappa_{3\mid2,1}}\) it is therefore
	sufficient to find a \(\magenta{\kappa_{3,2\mid 1}^g}\) that is only
	PIC w.r.t. \(g\) for each \(g\).
\end{remark}
\begin{proof}
	Let \(g\colon E_1\to D\) be a measurable function.
	By definition of the conditional expectation
	we need to show for all \(A_3\in \cE_3\) and all \(A_{1,2}\in \cE_1 \otimes
	\cE_2\)
	\[
		\E\Bigl[
			\blue{\ind_{A_{1,2}}(\xi_1, \xi_2^{(g(\xi_1))})}
			\teal{\kappa_{3\mid2,1}}\bigl(\xi_1, \xi_2^{(g(\xi_1))}, g(\xi_1); A_3\bigr)
		\Bigr]
		=\E\bigl[\blue{\ind_{A_{1,2}}(\xi_1, \xi_2^{(g(\xi_1))})}\ind_{A_3}(\xi_3)\bigr] 
	\]
	Without loss of generality we may only consider \(A_{1,2}= A_1\times A_2\in \cE_1\times \cE_2\)
	since the product sigma algebra \(\cE_1\otimes \cE_2\) is generated by these
	rectangles. 
	Since
	\[
		\kappa_{2\mid 1}^g(x_1; A_2)\coloneq \magenta{\kappa_{3,2\mid 1}}(x_1, g(x_1); E_3\times A_2)
	\]
	is a regular conditional version of \(\Pr(\xi_2^{(g(\xi_1))}\in \cdot \mid \xi_1)\)
	by assumption \eqref{eq: consistent outer kernel}
	we may apply disintegration \citep[e.g.][Thm. 6.4]{kallenbergFoundationsModernProbability2002} to the 
	measurable function
	\[
		\varphi(x_1, x_2) \mapsto \blue{\ind_{A_2}(x_2)}\teal{\kappa_{3\mid2,1}}(x_1, x_2, g(x_1); A_3)
	\]
	to obtain
	\begin{equation}
		\label{eq: disintegration of h}\begin{aligned}
		\E\bigl[\varphi\bigl(\xi_1, \xi_2^{g(\xi_1)}\bigr)\mid \xi_1\bigr]
		\overset{\as}&= \int \varphi(\xi_1, x_2)\kappa_{2\mid 1}^g(\xi_1; dx_2)
		\\
		\overset{\text{def.}}&= \int \varphi(\xi_1, x_2)\magenta{\kappa_{3,2\mid 1}}(\xi_1, g(\xi_1); E_3 \times dx_2).
	\end{aligned}	
	\end{equation}
	We thereby have
	\begin{align*}
		&\E\Bigl[
			\blue{\ind_{A_1}(\xi_1)\ind_{A_2}(\xi_2^{g(\xi_1)})}
			\teal{\kappa_{3\mid2,1}}\bigl(\xi_1, \xi_2^{g(\xi_1)}, g(\xi_1); A_3\bigr)
		\Bigr]
		\\
		&=\E\Bigl[\blue{\ind_{A_1}(\xi_1)}\varphi\bigl(\xi_1, \blue{\xi_2^{g(\xi_1)}}\bigr)\Bigr]
		\\
		\overset{\eqref{eq: disintegration of h}}&=
		\E\Bigl[
			\blue{\ind_{A_1}(\xi_1)}
			\int \varphi(\xi_1, \blue{x_2})
			\magenta{\kappa_{3,2\mid 1}}(\xi_1, g(\xi_1); E_3 \times dx_2)\Bigr]
		\\
		\overset{\text{Lemma}~\ref{lem: coupled kernel 3|2,1}}&=
		\E\Bigl[
			\blue{\ind_{A_1}(\xi_1)}
			\magenta{\kappa_{3,2\mid 1}}(\xi_1, g(\xi_1); A_3 \times \blue{A_2})
		\Bigr]
		\\
		\overset{\eqref{eq: consistent outer kernel}}&=
		\E\Bigl[
			\blue{\ind_{A_1}(\xi_1)}
			\ind_{A_3 \times \blue{A_2}}(\xi_3, \xi_2^{g(\xi_1)})
		\Bigr]
		\\
		&= \E\Bigl[
			\blue{\ind_{A_1}(\xi_1)\ind_{A_2}(\xi_2^{g(\xi_1)})}
			\ind_{A_3}(\xi_3)
		\Bigr].
	\end{align*}
	The crucial step is the application of Lemma~\ref{lem: coupled kernel
	3|2,1}, which provides an integral representation of a regular conditional
	distribution of \(\xi_3, \xi_2^y \mid \xi_1\) that \emph{couples} the two
	conditional kernels.

	\begin{lemma}\label{lem: coupled kernel 3|2,1}
		For all \(A_2\in \cE_2\), \(A_3\in \cE_3\), all \(y\in \domain\) and \(\Pr_{\xi_1}\)-almost all \(x_1\)
		\begin{align}
			\label{eq: coupled kernel 3|2,1}
			\magenta{\kappa_{3,2\mid 1}}(x_1, y;A_3 \times A_2)
			&=\int \varphi(x_1,x_2) \magenta{\kappa_{3,2\mid 1}}(x_1, y; E_3 \times dx_2)
			\\
			\nonumber
			&= \int \ind_{A_2}(x_2) \teal{\kappa_{3\mid 2,1}}(x_1, x_2, y; A_3) \magenta{\kappa_{3,2\mid 1}}(x_1, y; E_3 \times dx_2).
		\end{align}
	\end{lemma}
	In the remainder of the proof we will show this Lemma.
	To this end pick any \(A_1 \in \cE_1\). Then by definition of the
	conditional expectation \citep[e.g.\label{footnote: test}][chap.~8]{klenkeProbabilityTheoryComprehensive2014}
	\begin{align*}
		&\E[\blue{\ind_{A_1}(\xi_1)}\magenta{\kappa_{3,2\mid 1}}(x_1, y; A_3 \times A_2)]
		\\
		\overset{\eqref{eq: consistent outer kernel}}&= \E[\blue{\ind_{A_1}(\xi_1)} \ind_{A_2}(\xi_2^y) \ind_{A_3}(\xi_3)]
		\\
		\overset{\eqref{eq: joint inner kernel}}&=
		\E[\blue{\ind_{A_1}(\xi_1)} \ind_{A_2}(\xi_2^y) \teal{\kappa_{3\mid 2,1}}(\xi_1, \xi_2^y; A_3)]
		\\
		\overset{(*)}&= 
		\E\Bigl[\blue{\ind_{A_1}(\xi_1)} \int \ind_{A_2}(x_2) \teal{\kappa_{3\mid 2,1}}(\xi_1, x_2, y; A_3) \magenta{\kappa_{3,2\mid 1}}(\xi_1, y; E_3 \times dx_2)\Bigr]
	\end{align*}
	Note that the constant function \(g\equiv y\) is always measurable for the application of \eqref{eq: consistent outer kernel}.
	The last step \((*)\) is implied by disintegration \citep[e.g.][Thm. 6.4]{kallenbergFoundationsModernProbability2002}
	\[
		\E[f(\xi_1, \xi_2^y) \mid \xi_1] \overset{\as}= \int f(\xi_1, x_2) \kappa_{2\mid 1}^y(\xi_1; dx_2)
	\]
	of the measurable function
	\(
		f(x_1, x_2)
		\coloneq \ind_{A_2}(x_2)\teal{\kappa_{3\mid 2,1}}(x_1, x_2; A_3)
	\)
	using the probability kernel
	\[
		\kappa_{2\mid 1}^y(x_1; A_2)\coloneq \magenta{\kappa_{3,2\mid 1}}(x_1, y; E_3\times A_2)
	\]
	which is a regular conditional version of \(\Pr(\xi_2^y \in \cdot \mid \xi_1)\) by
	assumption \eqref{eq: consistent outer kernel}, since the constant function
	\(g\equiv y\) is measurable.
\end{proof}

\subsubsection{Proof of the main result (Theorem \ref{thm: conditionally independent sampling})}

In this section we assume that \(\domain\) is a locally compact,
separable and metrizable space and \(\range\) is a Polish space.
In this setting, we can use Proposition \ref{prop: dependent consistent}
to get rid of the assumption in Proposition \ref{prop: consistency shuffle}.
This is done in the following corollary, which is then modified for
conditional independence in Lemma \ref{lem: consistency allows conditional independence}.
We finally prove our main result using repeated applications of this lemma.

\begin{corollary}[Automatic PIC]
	\label{cor: automatic consistency}
	Let \(Z\) be a random variable in a standard borel space \((E, \borel(E))\),
	and \(\rf\) a random variable in \(C(\domain, \range)\).
	Let \(W\) be a random element in an arbitrary measurable space \((\Omega,
	\cF)\). If there exists a joint conditional distribution \(\kernel\) for
	\(Z\) given \((W, \rf(x))\) then
	\(\kernel\) is automatically PIC.
	That is, for all \(B\in \borel(E)\) all \(\sigma(W)\) measurable \(X\)
	\[
		\Pr(Z\in B \mid W, \rf(X)) \overset\as= \kernel(W, \rf(X), X; B),
	\]
\end{corollary}
\begin{proof}
		Since \(X=g(W)\) for some measurable function \(g\), Proposition
		\ref{prop: consistency shuffle} with \((\xi_3, \xi_2^y,
		\xi_1)=(Z,\rf(y), W)\) yields the claim, since a PIC joint
		probability kernel for \(\xi_3, \xi_2^y\) given \(\xi_1\) \emph{exists} by
		Proposition \ref{prop: dependent consistent}.
\end{proof}

The proof of Theorem \ref{thm: conditionally independent sampling} \ref{it: consistent without dependent variable} will follow from
repeated applications of the following lemma. This lemma slightly
generalizes Corollary \ref{cor: automatic consistency} to
the conditional independence setting.
Theorem \ref{thm: conditionally independent sampling} \ref{it: consistent with dependent variable}
will follow from \ref{it: consistent without dependent variable} and an application
of Proposition \ref{prop: dependent consistent}.

\begin{lemma}[PIC allows conditional independence]
	\label{lem: consistency allows conditional independence}
	Let \(Z\) be a random variable in a standard borel space \((E, \borel(E))\),
	\(\rf\) a continuous random function in \(C(\domain, \range)\). Let \(W\) be
	a random element in an arbitrary measurable space \((\Omega, \cF)\). If
	there exists a joint conditional distribution \(\kernel\) for \(Z\) given
	\(W, \rf(x)\) then for any random variable \(X \indep_W (Z,\rf)\) in
	\(\domain\)
	\[
		\Pr(Z\in B \mid W, X, \rf(X)) = \kernel(W, \rf(X), X; B).
	\]
\end{lemma}

\begin{proof}
	Observe that \(X\) is clearly measurable with respect to \(W^+ \coloneq (W,X)\). Our
	proof strategy therefore relies on constructing a joint conditional distribution
	for \(Z\) given \((W^+, \rf(x))\) using \(\kernel\) and apply Corollary \ref{cor: automatic
	consistency}.
	
	Since \(X\) independent from \((Z, \rf)\) conditional on \(W\)
	there exists a standard uniform \(U\sim \uniform(0,1)\) independent from \((W,Z,\rf)\) such
	that \(X = h(W,U)\) for some measurable function \(h\) \citep[Prop.\@
	6.13]{kallenbergFoundationsModernProbability2002}.
	Since \(U\) is independent from \(W,Z,\rf\) we have by \citep[Prop.\@ 6.6]{kallenbergFoundationsModernProbability2002}
	\[
		\Pr(Z\in B \mid W, U, \rf(x))
		= \Pr(Z\in B \mid W, \rf(x))
		= \kernel(W, \rf(x), x; B)
	\]
	for all \(x\in \domain\). Since \(\sigma(W,X,\rf(x)) \subseteq \sigma(W, U, \rf(x))\) and
	\((W,\rf(x))\) is measurable with respect to \(\sigma(W,X,\rf(x))\) we therefore have
	\[
		\Pr(Z\in B \mid \underbrace{W, X}_{=W^+}, \rf(x))
		= \kernel(W, \rf(x),x;B)
		\eqcolon \kernel^+(\underbrace{W, X}_{=W^+}, \rf(x),x;B),
	\]
	where \(\kernel^+\) is defined as constant in the second input.
	An application of Corollary \ref{cor: automatic consistency} to
	\(\kernel^+\) yields the claim.
\end{proof}

\begin{proof}[Proof of Theorem \ref{thm: conditionally independent sampling}]
	We will prove  \ref{it: consistent without dependent variable} by induction over \(k\in\{0,\dots,n+1\}\).
	For any conditionally independent evolution \((X_n)_{n\in \nat_0}\), the induction claim is
	\begin{equation}
		\label{eq: induction claim, conditionally indep. sampling}
		\begin{aligned}
		&\Pr\bigl(
			Z \in A
			\mid W, (\rf_i(X_i))_{i\in[0\inter k)},
			(\rf_i(x_i))_{i\in[k\inter n]}, X_{[0\inter k)}
		\bigr)
		\\
		&\qquad\quad= \kernel\bigl(W, (\rf_i(X_i))_{i\in [0\inter k)}, (\rf_i(x_i))_{i\in [k\inter n]}, X_{[0\inter k)}, x_{[k\inter n]}; A\bigr).
	\end{aligned}\end{equation}
	That is, we plugged in all the random variables up to index \(k-1\).
	The induction start with \(k=0\) is given by assumption and \(k=n+1\) is the claim,
	so we only need to show the induction step \(k\to k+1\). For this purpose we
	want to define \(\tilde W
	= \bigl(W, (\rf_i(X_i))_{i\in [0\inter k)}, (\rf_i(x_i))_{i\in(k\inter n]}, X_{[0\inter k)}\bigr)\)
	and the kernel
	\[
		\tilde\kernel_{x_{(k\inter n]}}\bigl(\tilde W, \rf_k(x_k), x_k; A\bigr)
		\coloneq \kernel\bigl(W, (\rf_i(X_i))_{i\in [0\inter k)}, (\rf_i(x_i))_{i\in [k\inter n]}, X_{[0\inter k)}, x_{[k\inter n]}; A\bigr),
	\]
	which is formally defined for any fixed \(x_{(k\inter n]}\) by mapping the elements of \(\tilde W\) into the right 
	position. By induction \eqref{eq: induction claim, conditionally indep. sampling} we thereby have
	\[
		\Pr(Z \in A \mid \tilde W, \rf(x_k))
		= \tilde\kernel_{x_{(k\inter n]}}(\tilde W, \rf(x_k), x_k; A).
	\]
	We can now finish the induction using Lemma \ref{lem: consistency allows conditional independence}
	if we can prove \(X_k\) is independent from \((Z, \rf)\) conditional on \(\tilde W\),
	because then we can also plug-in \(X_k\).
	For the conditional independence we will use the characterization in
	Proposition~6.13 by \citet{kallenbergFoundationsModernProbability2002}.

	Since \(X_k\) is independent from \((Z,\rf)\) conditional on \(\filt_{k-1}\) there exists,
	by this Proposition, a uniform random variable \(U\sim \uniform(0,1)\)
	independent from \((Z,\rf, \filt_{k-1})\)
	such that \(X_k = h(\xi, U)\) for some measurable function \(h\) and a random
	element \(\xi\) that generates \(\filt_{k-1}\). Due to \(\filt_{k-1}\subseteq \sigma(\tilde W)\)
	the element \(\xi\) is a measurable function of \(\tilde W\) and therefore \(X_k=
	\tilde h(\tilde W, U)\) for some measurable function \(\tilde h\).
	Since \(U\) is independent from \((Z,\rf,\filt_{k-1})\), it is independent from
	\(\tilde W\) as \(\sigma(\tilde W)\subseteq \sigma(\rf, \filt_{k-1})\). Using Prop.~6.13 from \citet{kallenbergFoundationsModernProbability2002}
	again, \(X_k\) is thereby independent from \((Z, \rf)\) conditional on \(\tilde W\).

	What remains is the proof of \ref{it: consistent with dependent variable}.
    Let \(x_n\) be fixed and define \(\tilde Z=(Z,\rf_n(x_n))\). Since \(\kernel\) is a joint conditional
    distribution for \(Z,\rf_n(x_n)\) given \(\filt, (\rf_k(x_k))_{k\in [0\inter
    n)}\) the kernel
    \[
        \kernel_{x_n}(y_{[0\inter n)}, x_{[0\inter n)}; B) \coloneq  \kernel(y_{[0\inter n)}, x_{[0\inter n]}; B)
    \]
    clearly satisfies the requirements of \ref{it: consistent without dependent variable} and thereby
    \[
        \Pr\bigl(Z,\rf(x_n) \in B \mid \filt_{n-1}\bigr)
        = \kernel\bigl(\underbrace{W, (\rf_k(X_k))_{k\in [0\inter n)}, X_{[0\inter n)}}_{\eqcolon\tilde W}, x_n; B\bigr)
        \eqcolon \tilde\kernel(\tilde W, x_n; B).
    \]
    Since \(\tilde W\) generates \(\filt_{n-1}\) we are almost in the setting of
    Proposition \ref{prop: dependent consistent}, as we have
    continuity in \(x_n\). However, since \(X_n\) is not previsible we have to
    repeat the same trick we used in the proof of Lemma \ref{lem: consistency
    allows conditional independence}. Namely, \(X_n\) is measurable with
    respect to \(W^+ \coloneq (\tilde W, X_n)\) and we will have to construct a
    joint conditional distribution for \(Z,\rf(x_n)\) given \(W^+\).

    Since \(X_n\) is independent from \(Z,\rf\) conditional on \(\tilde W\),
    there exists a standard uniform \(U\sim \uniform(0,1)\) independent from
    \((\tilde W, Z,\rf)\) such that \(X_n = h(\tilde W, U)\) for some measurable
    function \(h\) \citep[Prop.\@
    6.13]{kallenbergFoundationsModernProbability2002}.  Since \(U\) is
    independent from \((\tilde W, Z, \rf)\), we have by
    \citep[Prop.\@ 6.6]{kallenbergFoundationsModernProbability2002}
    \[
        \Pr(Z,\rf(x_n) \in B \mid \tilde W, U)
        \overset\as= \Pr(Z,\rf(x_n) \in B \mid \tilde W)
        \overset\as= \tilde\kernel(\tilde W, x_n; B)
    \]
    for all \(x_n\in \domain\). Since \(\sigma(\tilde W,X_n) \subseteq \sigma(\tilde W, U)\) and
	\(\tilde W\) is measurable with respect to \(\sigma(\tilde W,X_n)\) we therefore have
	\[
		\Pr(Z,\rf(x_n)\in B \mid \underbrace{\tilde W, X_n}_{=W^+})
		\overset\as= \tilde\kernel(\tilde W,x;B)
		\eqcolon \kernel^+(\underbrace{\tilde W, X_n}_{=W^+}, x_n;B),
	\]
	where \(\kernel^+\) is defined as constant in the second input. Clearly, by
    definition of \(\kernel^+\) via \(\kernel\), \(\kernel^+\) is continuous in \(x_n\)
    and as a continuous joint conditional distribution it is PIC by
    Proposition \ref{prop: dependent consistent}.
    This finally implies the claim
    \begin{align*}
        \Pr(Z,\rf(X_n) \in B \mid \smash{\underbrace{\filt_{n-1}, X_n}_{=W^+}})
        \overset\as&= \kernel^+(W^+, X_n; B)
        \\
        &= \kernel(W, (\rf_k(X_k))_{k\in [0\inter n)}, X_{[0\inter n)}, X_n; B).
        \qedhere
    \end{align*}
\end{proof}

	\section*{Acknowledgements}
	I want to thank my colleagues at the University of Mannheim -- especially my
	PhD Supervisor Leif Döring and Martin Slowik -- for insightful discussions and
	support.

	\bibliography{references}

	\appendix
\section{Gaussian conditionals}
\label{sec: gaussian conditionals}

In this section we want to recall the canonical version of the conditional
distribution of Gaussian random vectors/functions.

\begin{remark}[Multivariate output]
    While we limit ourselves to real-valued functions in the following lemma,
    the result may easily be extended to \(\real^\dims\) valued functions. To do
    so simply turn the output dimension into an input, i.e.\ \(\rf(x)_i =:
    \rf(x,i)\) for \(i \in [1\inter \dims]\), where we recall the notation for
    discrete intervals
    \[
[i\inter j] \coloneq [i,j] \cap \integer,
        \qquad
        [i\inter j) \coloneq [i,j)\cap \integer,
        \qquad
        \text{etc.}
    \]
    In other words we simply replace the domain
    \(\domain\) by \(\domain\times [1\inter \dims]\).
    This necessitates that we do not only consider \(\rf(x_n)\)
    given \(\rf(x_0), \dots, \rf(x_{n-1})\) but instead the the distribution
    of \(\rf(x_k), \dots, \rf(x_n)\) given \(\rf(x_0), \dots, \rf(x_{k-1})\).
\end{remark}

\begin{lemma}[Gaussian conditional distribution]
    \label{lem: gaussian conditional}
    Let \((\rf(x))_{x\in \domain}\) with values in \(\real\) be a Gaussian random function
    with mean and covariance function
    \[
        \mu(x) = \E[\rf(x)]
        \qquad\text{and}\qquad
        \C_\rf(x,y) = \Cov(\rf(x), \rf(y)).
    \]
    Then a regular conditional distribution
    \[
        \Pr\Bigl(\bigl(\rf(x_{k+1}), \dots, \rf(x_n)\bigr)\in A \bigm| Y\Bigr)
        = \kernel_x(Y; A)
    \]
    for \(k< n\), \(Y=(\rf(x_0), \dots, \rf(x_k))\) and \(x=(x_0, \dots, x_n)\) is given by the Gaussian distribution
    \begin{align}
        \kernel_{x}(y; A)
        &= \kernel(x, y; A)
        \nonumber
        \\
        \label{eq: canonical gaussian conditional}
        &\propto \int_A \exp\Bigl(-\tfrac12\bigl(t-\mu_{n\mid k}(x,y)\bigr)^\transpose [\Sigma_{n\mid k}(x)]^{-1}\bigl(t-\mu_{n\mid k}(x,y)\bigr)\Bigr)dt,
    \end{align}
    where the posterior mean is given by
    \[
        \mu_{n\mid k}(x, y)
        = \begin{pmatrix}
            \mu(x_{k+1})
            \\
            \vdots
            \\
            \mu(x_n)
        \end{pmatrix} +
        M \Bigl[(\C_\rf(x_i, x_j))_{i,j\in [0\inter k]}\Bigr]^{-1}
        \Bigl(y - \begin{pmatrix}
            \mu(x_0)\\
            \vdots\\
            \mu(x_k)
        \end{pmatrix}\Bigr)
    \]
    and posterior covariance by
    \[
        \Sigma_{n\mid k}(x)
        = (\C_\rf(x_i, x_j))_{i,j\in (k\inter n]}
        - M \Bigl[(\C_\rf(x_i, x_j))_{i,j\in [0\inter k]}\Bigr]^{-1} M^\transpose
    \]
    with
    \[
        M = \begin{pmatrix}
            \C_\rf(x_{k+1}, x_0) & \cdots & \C_\rf(x_{k+1}, x_k)\\
            \vdots & & \vdots\\
            \C_\rf(x_n, x_0) & \cdots & \C_\rf(x_n, x_k)
        \end{pmatrix}.
    \]
\end{lemma}
\begin{proof}
    Since \((\rf(x_0), \dots, \rf(x_n))\) is multivariate Gaussian,
    this is simply a standard result about conditionals of Gaussian vectors
    \citep[e.g.][Prop.\ 3.13]{eatonMultivariateStatisticsVector2007}.
\end{proof}

\begin{definition}[Canonical Gaussian conditional distribution]
    \label{def: canonical Gaussian dist}
    Since the definition of the regular conditional distribution in \eqref{eq: canonical gaussian conditional}
    is not unique, we refer to this specific version as the \emph{canonical} Gaussian conditional distribution.
\end{definition}

\begin{remark}[Joint kernel]
    \label{rem: Gaussian join kernel}
    The collection of kernels \((\kernel_x)_{x\in \domain^{n+1}}\) is a
    joint distribution, since it is measurable in \(x,y\) for fixed
    \(A\).
\end{remark}

\begin{remark}[Continuity]
    \label{rem: continuity of gaussian conditional}
    If \(\mu\) and \(\C_\rf\) are continuous, then the characteristic function
    of the joint kernel
    \[
        \hat{\kernel}(x,y; u) = \exp\Bigl(i u^\transpose \mu_{n\mid k}(x,y) - \tfrac12 u^\transpose \Sigma_{n\mid k}(x) u\Bigr)
    \]
    is continuous in \(x_{(k\inter n]}\), since both \(\mu_{n\mid k}\) and \(\Sigma_{n\mid k}\) are
    continuous in \(x_{(k\inter n]}\). This can be seen directly from the explicit
    expressions in Lemma \ref{lem: gaussian conditional} as the \(x_j\) with
    \(j > k\) are not involved in the matrix inversion. This implies
    continuity of \(x_{(k\inter n]} \mapsto \kernel(x,y; \cdot)\) in the weak 
    topology by Lévy's continuity theorem \citep[e.g.][Thm.\
    5.3]{kallenbergFoundationsModernProbability2002}.
\end{remark}

\begin{remark}[Artificially breaking things]
    This canonical kernel may be artificially modified to be
    discontinuous. E.g.\ with an indicator on \(\{x=y\}\), which is a null set for
    every fixed \(x\) as \(Y\) has a density. The resulting collection of kernels
    would still be a joint regular conditional distribution, however
    in contrast to the canonical Gaussian conditional distribution
    they would not necessarily be PIC (see Prop.~\ref{prop: dependent consistent}).
    Similarly, the collection may be modified to not be a joint kernel.
\end{remark}
 	\section{Topological foundation}
\label{sec: topological foundation}

In this section we show the evaluation function to be continuous and therefore
measurable for continuous random functions.
For \emph{compact} \(\domain\) this result can be collected from various
sources \citetext{e.g.\@ \citealp[Thm.\@ 4.2.17]{engelkingGeneralTopologyRevised1989}
and \citealp[Thm.\@ 4.19]{kechrisClassicalDescriptiveSet1995}}. But we could not
find a reference for the result in this generality, so we provide a proof.

\begin{theorem}[Continuous functions]
    \label{thm: topology of continuous functions}
    Let \(\domain\) be a locally compact, separable and metrizable space\footnote{\label{footnote: regular second countable}We do not need \(\domain\) to be metrizable but only regular
        and second countable (in metrizable spaces `second countable' is
        equivalent to `separable' \citep[Cor.\@
        4.1.16]{engelkingGeneralTopologyRevised1989}). We choose this more
        specific definition to make it more obvious that a locally compact
        Polish
        space satisfies the requirements. However, for the proof we will assume
        the more general setting.
    },
    \(\range\) a Polish space and \(C(\domain, \range)\) the space of continuous
    functions equipped with the \emph{compact-open}\footnote{
        The sets \(M(K,U)\coloneq \{f\in C(\domain, \range): f(K) \subseteq U\}\) with
        \(K\subseteq \domain\) compact and \(U\subseteq \range\) open, form a
        sub-base of the \emph{compact-open} topology \citep[e.g.][Sec.
        3.4]{engelkingGeneralTopologyRevised1989}. I.e.\@ the compact-open
        topology it is the smallest topology such that all \(M(K,U)\) are open.
        Recall that the set of finite intersections of a sub-base form a base of the topology
        and elements from the topology can be expressed as unions of base elements. 
    } topology. Then
    \begin{enumerate}[label=\normalfont{(\roman*)}]
        \item\label{it: eval function is continuous}
        the evaluation function
        \[
            e\colon \begin{cases}
                C(\domain, \range) \times \domain \to \range
                \\
                (f,x) \mapsto f(x)
            \end{cases}
        \]
        is continuous and therefore measurable.
        \item\label{it: continuous function polish with metric}
        \(C(\domain,\range)\) is a \textbf{Polish space}, whose topology is
        generated by the metric
        \[
            \metric(f,g) \coloneq \sum_{n=1}^\infty 2^{-n} \frac{\metric_n(f,g)}{1+\metric_n(f,g)}
            \quad \text{with}\quad
            \metric_n(f,g) \coloneq \sup_{x\in K_n} \metric_\range(f(x), g(x))
        \]
        for any metric \(\metric_\range\) that generates the topology of \(\range\)
        and any \emph{compact exhaustion}\footnote{
            The set \(\domain\) is \emph{hemicompact} if it can be
            \emph{exhausted by the compact sets} \((K_n)_{n\in \nat}\), which
            means that the compact set \(K_n\) is contained in the interior of
            \(K_{n+1}\) for any \(n\) and \(\domain =\bigcup_{n\in \nat} K_n\).
        } \((K_n)_{n\in \nat}\) of \(\domain\),
        that always exists because \(\domain\) is hemicompact!
        \item\label{it: borel coincides} The Borel \(\sigma\)-algebra of \(C(\domain, \range)\) is equal to
        the restriction of the product sigma algebra of \(\range^\domain\) to
        \(C(\domain, \range)\), i.e. \(\borel(C(\domain, \range)) =
        \borel(\range)^{\otimes \domain}\bigr|_{C(\domain, \range)}\).
    \end{enumerate}
\end{theorem}

\begin{remark}[Topology of pointwise convergence]
    \label{rem: top of pointwise convergence}
    The topology of point-wise convergence ensures that all
    projection mappings \(\pi_x(f) = f(x)\) are continuous. It coincides with
    the product topology \citep[Prop.
    2.6.3]{engelkingGeneralTopologyRevised1989}.
    Thm.~\ref{thm: topology of continuous functions} \ref{it: borel coincides} ensures
    that the Borel-\(\sigma\)-algebra generated by the topology of
    point-wise convergence coincides with the Borel
    \(\sigma\)-algebra generated by the compact-open
    topology.
\end{remark}
\begin{remark}[Construction]
    The main tool for the construction of probability measures, Kolmogorov's
    extension theorem \citep[e.g.][Sec.\@ 14.3]{klenkeProbabilityTheoryComprehensive2014}, allows for the construction of random measures
    on product spaces. This is only compatible with the product topology, i.e.\
    the topology of point-wise convergence. But the evaluation map is generally
    not continuous with respect to this topology \citep[Prop.\@
    2.6.11]{engelkingGeneralTopologyRevised1989}. \ref{it: borel coincides}
    ensures that this does not pose a problem as long as \(\domain\) and
    \(\range\) satisfy the requirements of Theorem~\ref{thm: topology of
    continuous functions} and the constructed random process has a continuous
    version \citetext{cf.\@ \citealp[Thm.
    3]{talagrandRegularityGaussianProcesses1987},
    \citealp{costaSamplePathRegularity2024} and references therein}.
\end{remark}

\begin{remark}[Limitations]
    \label{rem: limits of the result}
    While the compact-open topology can be defined for general topological
    spaces, the continuity of the evaluation map crucially depends on \(\domain\)
    being locally compact \citep[Thm.\@ 3.4.3 and comments
    below]{engelkingGeneralTopologyRevised1989}. For \(\domain\) and \(\range\)
    Polish spaces, this implies \(C(\domain,\range)\) is generally only 
    well behaved if \(\domain\) is locally compact.
\end{remark}

\begin{remark}[Discontinuous case]
    \label{rem: discontinuous case}
    Without continuity it is already difficult to obtain a random
    function \(\rf\) that is almost surely measurable and can be evaluated point-wise. The construction of
    Lévy processes in càdlàg\footnote{
        french: continue à droite, limite à gauche, ``right-continuous with left-limits''
    } space only works on ordered domains such as \(\real\), where
    `right-continuous' has meaning. Typically, discontinuous random
    functions are therefore only constructed as generalized functions in the
    sense of distributions\footnote{
        The set of distributions is defined as the topological dual to a
        set of test functions. In particular, distributions are \emph{continuous}
        linear functionals acting on the test functions. Thereby one may
        hope that Theorem~\ref{thm: topology of continuous functions} is applicable, but
        the set of test functions is typically not locally compact (cf.\@ Remark \ref{rem: limits of the result}).
    } that cannot be evaluated point-wise
    \citep[e.g.][]{schafflerGeneralizedStochasticProcesses2018}. In particular, we
    cannot hope to evaluate generalized random functions at random locations.
    Nevertheless it may be possible to prove the evaluation function to be
    measurable for more general separable functions.
\end{remark}

\begin{proof}
    Since \(\domain\) is locally compact, \ref{it: eval function is continuous}
    follows from Proposition 2.6.11 and Theorem 3.4.3. by
    \citet{engelkingGeneralTopologyRevised1989}.
    
    For \ref{it: continuous function polish with metric} let us begin to show that
    \(\domain\) is \textbf{hemicompact/exhaustible by compact sets}. Since the space \(\domain\) is
    locally compact, pick a compact neighborhood for every point. The interiors
    of these compact neighborhoods obviously cover \(\domain\). Since every
    regular, second countable space\footref{footnote: regular second countable} is Lindelöf \citep[Thm.\@ 3.8.1]{engelkingGeneralTopologyRevised1989}, we can pick a countable
    subcover, such that the interiors of the sequence \((C_i)_{i\in \nat}\) of
    compact sets cover the domain \(\domain\). We inductively
    define a compact exhaustion \((K_n)_{n\in \nat}\) with \(K_1 \coloneq C_1\).
    Observe that the set \(K_n\) is covered by the interiors \((\interior{C_i})_{i\in
    \nat}\). Since \(K_n\) is compact, we can choose a finite sub-cover
    \((\interior{C_i})_{i\in I}\) and define \(K_{n+1} \coloneq \bigcup_{i \in I} C_i\cup C_{n+1}\).
    Then by definition \(K_n\) is contained in the interior of the compact set
    \(K_{n+1}\) and due to \(C_n \subseteq K_n\) this sequence also covers the space \(\domain\)
    and is thereby a compact exhaustion.

    It is straightforward to check that the metric defined in \ref{it:
    continuous function polish with metric} is a metric, so we will only prove
    this \textbf{metric induces the compact-open} topology.
    \begin{enumerate}[wide,label={(\Roman*)}]
        \item \textbf{The compact-open topology is a subset of the metric topology.}
        \label{it: compact-open is subset of metric top}
        We need to show that the sets \(M(K,U)\) are open with respect to the metric. 
        This requires for any \(f\in M(K,U)\) an \(\epsilon>0\) such that the
        epsilon ball \(B_\epsilon(f)\) is contained in \(M(K,U)\).

        We start by constructing a finite cover of \(f(K)\).
        For any \(x\in K\) there exists \(\delta_x>0\) with 
        \(B_{2\delta_x}(f(x))\subseteq U\) for balls induced by the metric
        \(\metric_\range\) as \(U\) is open. Since \(K\) is compact,
        \(f(K)\subseteq U\) is a compact set covered by the balls
        \(B_{\delta_x}(f(x))\). This yields a finite subcover
        \(B_{\delta_1}(f(x_1)), \dots, B_{\delta_m}(f(x_m))\) of \(f(K)\).

        Using this cover we will prove the following criterion: Any \(g\in C(\domain, \range)\) is in \(M(K, U)\) if
        \begin{equation}
            \label{eq: sup bound}
            \sup_{x\in K} \metric_\range(f(x), g(x))
            < \delta
            \coloneq \min\{\delta_1,\dots, \delta_m\}.
        \end{equation}
        For this criterion note that for any \(x\in K\) there exists \(i\in \{1,\dots, m\}\) such that \(f(x) \in B_{\delta_i}(f(x_i))\).
        This implies 
        \[
            \metric(g(x), f(x_i)) \le \metric(g(x), f(x)) + \metric(f(x), f(x_i)) \le 2\delta_i,
        \]
        which implies \(g(K)\subseteq \bigcup_{i=1}^m B_{2\delta_i}(f(x_i))\subseteq
        U\) and therefore \(g\in M(K,U)\).

        Consequently, if there exists \(\epsilon>0\) such that \(g\in B_\epsilon(f)\) implies
        criterion \eqref{eq: sup bound}, then we have \(B_\epsilon(f)
        \subseteq M(K,U)\) which finishes the proof. And this is what we will show.
        Since \(K\) is compact and the interiors of \(K_n\) cover the space, there exists
        a finite sub-cover \(K\subseteq \bigcup_{i\in I}K_i\) and therefore some \(m=\max I\) such that \(K\) is in the
        interior of \(K_m\). By definition of \(\metric_m\) it is thus clearly sufficient
        to ensure \(\metric_m(f,g) < \delta\). And since \(\varphi(x)= \frac{x}{1+x}\) is a strict
        monotonous function \(\epsilon \coloneq 2^{-m} \varphi(\delta)\) does the job, since
        \(
            2^{-m}\varphi(\metric_m(f,g))\le \metric(f,g) \le \epsilon
        \)
        implies \(\metric_m(f,g) \le \delta\).

        \item \textbf{The metric topology is a subset of the compact-open topology.}
        Since the balls \(B_\epsilon(f)\) form a base of the metric topology it is sufficient
        to prove them open in the compact-open topology. If 
        for any \(g\in B_\epsilon(f)\) there exists a compact \(C_1,\dots, C_m\subseteq \domain\)
        and open \(U_1,\dots, U_m\subseteq \range\) such that \(g\in \bigcap_{j=1}^m M(C_j, U_j) \subseteq B_\epsilon(f)\),
        then the ball is open since these finite intersections are open sets in
        the compact-open topology and their union over \(g\) remains open. But
        since there exists \(r>0\) such that \(B_r(g)\subseteq B_\epsilon(f)\),
        it is sufficient to prove for any \(r>0\) that there exist compact
        \(C_j\) and open \(U_j\) such that
        \begin{equation}
            \label{eq: compact open subset}
            \textstyle
            g\in V\coloneq\bigcap_{j=1}^m M(C_j, U_j) \subseteq B_r(g).
        \end{equation}
        For this purpose pick \(K_N\) from the compact exhaustion with sufficiently large \(N\) such that
        \(2^{-N} < \frac{r}2\).
        Pick a finite cover \(O_{x_1}, \dots O_{x_m}\) of \(K_N\) from the cover
        \(\{O_x\}_{x\in K_N}\) with \(O_x \coloneq g^{-1}(B_{r/5}(g(x)))\)
        and define the sets
        \[
            C_j \coloneq \closure{O_{x_j}} \cap K_N\qquad U_j \coloneq B_{r/4}(f(x_j)).
        \]
        Clearly the \(U_j\) are open and the \(C_j\) are compact and we will now
        prove they satisfy \eqref{eq: compact open subset}. Observe that \(g\in V\)
        since for all \(j\)
        \[
            g(C_j)\subseteq g(\closure{O_{x_j}}) \subseteq \closure{B_{r/5}(f(x_j))}\subseteq U_j.
        \]
        Pick any other \(h\in V\). Then for all \(x\in K_N\) there exists \(i\)
        such that \(x\in O_{x_i}\subseteq C_i\) and by definition of \(V\) this implies
        \(h(x)\in U_i\) and also \(g(x)\in U_i\) and thereby
        \(\metric_\range(h(x), g(x)) \le r/2\). This uniform bound implies \(\metric_N(h,g) \le r/2\) and therefore
        \[
            \metric(g,h) \le \Bigl(\sum_{n=1}^N 2^{-n} \metric_n(g,h)\Bigr) + 
            \Bigl(\sum_{n=N+1}^\infty 2^{-n}\Bigr)
            \le \metric_N(g,h) + 2^{-N} < r,
        \]
        since \(\metric_n(f,g)\le \metric_N(f,g)\) for \(n\le N\).
        Thus \(h\in B_r(g)\) which proves \eqref{eq: compact open subset}.
    \end{enumerate}
    As \(C(\domain, \range)\) is clearly metrizable, what is left to prove are its
    separability and completeness. Separability could be proven directly similarly to the
    proof of Theorem 4.19 in \citet{kechrisClassicalDescriptiveSet1995} but for
    the sake of brevity this result follows from Theorem 3.4.16 and Theorem
    4.1.15 (vii) by \citet{engelkingGeneralTopologyRevised1989} and the fact that
    \(\domain\) and \(\range\) are second countable. Completeness follows from the fact that any
    Cauchy sequence \(f_n\) induces a Cauchy sequence \(f_n(x)\) for any \(x\) by
    definition of the metric. And by completeness of \(\range\) there must exist a
    limiting value \(f(x)\) for any \(x\). The continuity of \(f\) follows from
    the uniform convergence on compact sets, since every compact set is contained
    in some \(K_n\) from the compact exhaustion (cf. last paragraph in \ref{it:
    compact-open is subset of metric top}).
    
    What is left to prove is \ref{it: borel coincides}. Since the projections are continuous with respect to the compact open
    topology, they are measurable with respect to the Borel-\(\sigma\)-algebra.
    The product sigma algebra, which is the smallest sigma algebra to ensure all
    projections are measuralbe, restricted to the continuous functions is
    therefore a subset of the Borel \(\sigma\)-algebra. To prove the opposite inclusion,
    we need to show that the open sets are contained in the product \(\sigma\)-algebra.
    Since the space is second countable \citep[Cor.\@
    4.1.16]{engelkingGeneralTopologyRevised1989} and every open set thereby a
    countable union of its base, it is sufficient to check that the open ball
    \(B_\epsilon(f_0)\) for \(\epsilon>0\) and \(f_0\in C(\domain, \range)\) is in
    the product sigma algebra restricted to \(C(\domain, \range)\). But since
    \(B_\epsilon(f_0) = H^{-1}([0,\epsilon))\) with \(H(f)\coloneq\metric(f,f_0)\), it is
    sufficient to prove \(H\) is \(\sigma(\pi_x : x\in \domain)\)-\(\borel(\real)\)-measurable,
    where \(\pi_x\) are the projections. \(H\) is measurable if \(H_n(f) =
    \metric_n(f,f_0)\) is measurable, as a limit, sum, etc.\@ \citep[Thm.\@ 1.88-1.92]{klenkeProbabilityTheoryComprehensive2014} of measurable functions.
    But since \(\domain\) is separable \citep[Cor.\@ 1.3.8]{engelkingGeneralTopologyRevised1989}, i.e.\@ has a countable dense subset \(Q\),
    we have by continuity of \(f\) and \(f_0\)
    \[
        \metric_n(f, f_0) = \sup_{x\in K_n} \metric_\range(f(x), f_0(x))
        = \sup_{x\in K_n\cap Q} \metric_\range(\pi_x(f), \pi_x(f_0)).
    \]
    Since \(\metric_\range\) is continuous and thereby measurable \citep[Thm.\@
    1.88]{klenkeProbabilityTheoryComprehensive2014},
    \(H_n\) is measurable as a countable supremum of measurable
    functions \citep[Thm.\@ 1.92]{klenkeProbabilityTheoryComprehensive2014}.
\end{proof}

\end{document}